\documentstyle[amssymb,amsfonts,12pt]{amsart}

\newtheorem{theorem}{Theorem}[section]
\newtheorem{proposition}[theorem]{Proposition}
\newtheorem{lemma}[theorem]{Lemma}
\newtheorem{corollary}[theorem]{Corollary}

\newtheorem{conjecture}[theorem]{Conjecture}

\def\C{{\mbox{\rm\kern.24em
\vrule width.03em height1.43ex depth-.052ex \kern-.26em C}}}
\def\QSet{\mbox{\rm\kern.24em
\vrule width.03em height1.48ex depth-.051ex \kern-.26em Q}}
\def\Z{{\bf Z}}
\def\R{{\mbox{\rm I\kern-.22em R}}}
\def\arcsinh{{\rm arcsinh}}

\def\P{{\bf P}}

\def\T{{\bf T}}

\def\size{{\rm size}}

\def\k{{\bf k}}

\def\p{{\bf p}}
\def\q{{\bf q}}

\def\id{{\rm id}}

\def\be#1{\begin{equation}\label{#1}}
\def\bas{\begin{align*}}
\def\eas{\end{align*}}
\def\bi{\begin{itemize}}
\def\ei{\end{itemize}}
\newenvironment{proof}{\noindent {\bf Proof} }{\endprf\par}
\def \endprf{\hfill  {\vrule height6pt width6pt depth0pt}\medskip}
\def\emph#1{{\it #1}}

\title[A Carleson type theorem]{A Carleson type theorem for a Cantor group model of the
scattering transform}

\author{Camil Muscalu}
\address{Department of Mathematics, UCLA, Los Angeles CA 90095-1555}
\email{camil@@math.ucla.edu}

\author{Terence Tao}
\address{Department of Mathematics, UCLA, Los Angeles CA 90095-1555}
\email{tao@@math.ucla.edu}

\author{Christoph Thiele}
\address{Department of Mathematics, UCLA, Los Angeles CA 90095-1555}
\email{thiele@@math.ucla.edu}

\begin{document}

\begin{abstract}  

We consider a basic $d$-adic model for the scattering transform
on the line. We prove 
$L^2$ bounds for this scattering transform
and a weak $L^2$ bound for a Carleson type maximal operator
(Theorem \ref{main}).
The latter implies boundedness of $d$-adic models 
of generalized eigenfunctions of Dirac type operators with 
potential in 
$L^2(\R)$. We show that 
this result cannot be obtained by estimating the
terms in the natural multilinear expansion of the scattering
transform (Proposition \ref{unbounded-3}).

\end{abstract}

\maketitle

\section{Introduction}

It is widely understood that scattering transforms
are non-linear variants of the one dimensional Fourier transform. 
Thus scattering transforms give nonlinear Fourier transforms of
scalar or more generally matrix valued potentials $F(x)$.
For harmonic analysts this suggests to study the basic a priori estimates 
in Fourier analysis (such as for example Hausdorff Young 
inequalitites or estimates for
the Carleson operator) in the case of 
the scattering transforms. This naturally leads to the study of
the nonlinear Fourier transform for rough and slowly decaying
potentials. Beals and Coifman \cite{bealscoifman}
study in detail the case when the potential is (generic)
in $L^1$ or in weighted spaces $L^1\cap L^2$ with weights
of the form $(1+|x|)^m$ for suitable $m$.
More recently, Christ and Kiselev 
\cite{ck0}, \cite{ck1} have proven analogues of the Hausdorff Young inequality and a maximal Hausdorff Young inequality for a scattering transform. 
This is an estimate for potentials in $L^p$. Their result implies 
boundedness of 
eigenfunctions of one dimensional Schr\"odinger operators with potential in
$L^p$, $1<p<2$ for almost all positive energies. By an extension by 
Simon of a theorem of Sch'nol \cite{schnol} this implies that
the absolutely continuous spectrum of the Schr\"odinger operator
is supported on the entire positive half axis, see also 
\cite{simon} page 501. This implication
was one of the motivations of Christ and Kiselev to study
the maximal Hausdorff Young inequalities for the scattering transform.
We propose to study the analogue of Carleson's theorem \cite{carleson} 
or the sharper form by Hunt \cite{hunt}, see also \cite{lt} for a 
recent proof, for scattering transforms. This amounts to an $L^2$ endpoint of
the results by Christ and Kiselev and would give boundedness of
eigenfunctions of Schr\"odinger operators with potential in $L^2$.
The question of absolutely continuous spectrum for potentials in $L^2$
has been settled to the affirmative by Deift and Killip \cite{deiftkillip}, 
but this is a weaker statement than the conjectured boundedness
of eigenfunctions.

Currently we are not able to decide
whether the analogue of Carleson's theorem as stated below
is true or false.
The purpose of this article is to study a $d$-adic model
for this problem and prove a positive result for this model.

We restrict attention to one of the easiest cases of
the scattering transform. Thus consider the special AKNS-ZS system
(named after \cite{akns} and \cite{zakharovshabat}):
\begin{equation}\label{akns}
\frac{d f}{d x}= k J f(x) + q(x) f(x)
\end{equation}
for the unknown function $f:\R\to \C^2$ where $k\in \C$
is a spectral parameter,
$$J=\left(\begin{array}{cc} -i & 0 \\ 0 & i \end{array}\right)\ ,\ \ 
q(x)=\left(\begin{array}{cc} 0 & {F(x)} \\  \overline{F(x)} & 0 \end{array}\right)\ \ .$$ 
This can be read as eigenfunction equation for Dirac operators
on the real line. More generally one can write the eigenfunction equation
for Schr\"odinger operators on the real line in the framework of
AKNS-ZS systems. This links to the work of Christ and Kiselev,
but we shall not elaborate on this generalization. We remark that
Conjecture \ref{weak-carleson-conjecture} (as well as the
other conjectures formulated below) would imply boundedness
of solutions to (\ref{akns}) for almost every $k\in \R$.

We shall assume that $F$ is locally integrable
and for simplicity compactly supported.
Writing $a(x)\exp(-ikx)$ and $b(x)\exp(ikx)$ for the two components
of $f$ and assuming $k$ is real we obtain
the following equivalent ordinary differential equation
\begin{equation}\label{equivalent-ode}
G'=WG
\end{equation}
where
\begin{equation}\label{nonlinfour}
G=\left(\begin{array}{cc} a &
{b}
\\ \overline{b} &
\overline{a}
\end{array}\right)\ ,\ \ 
W(x)=\left(\begin{array}{cc} 0 & F(x)\exp(2ikx)
\\  \overline{F}(x)\exp(-2ikx) & 0
\end{array}\right)\ \ .
\end{equation}

Since $F$ is compactly supported, equation (\ref{equivalent-ode})
forces $G$ to be constant near
$-\infty$ and near $\infty$. Let $G(-\infty)$ and $G(\infty)$
denote these constant values.
Imposing the initial condition $G(-\infty)=\id$, then
standard existence theorems give a unique absolutely continuous solution
satisfying (\ref{equivalent-ode}) almost everywhere. Thus we can 
define $G(\infty)$ to  
be the scattering transform of
the potential $F$ at the spectral value $k$.

We have implicitly used that the differential equation (\ref{equivalent-ode})
forces $G$ to remain in the form stated in (\ref{nonlinfour}),
if it is initially of that form. It also forces $G$ to have
constant determinant, which for our chosen initial condition is 
equal to $1$. In other words, $G$ takes values in the Lie group 
$SU(1,1)$, see also the discussion in Section 
\ref{plancherel-section}.

To prove a priori estimates for the scattering transform,
we need a notion of size for the matrices $G$. A natural
size appearing in the $L^2$ theory of this scattering transform is
$\sqrt{\log|a|}$ where $a$ is the upper left entry of the 
matrix $G$. Observe that this quantity is positive since
$|a|^2-|b|^2=1$.

The
following are known analogues of standard estimates for
the Fourier transform:
Recall that $a(\infty)$ and $a(x)$ for given $x$ are functions
in the parameter $k$ which we have suppressed in the notation.

\begin{theorem}\label{classical}
Riemann - Lebesgue estimate
$$\|    \sqrt{\log |a(\infty)|} \|_\infty \le C \|F\|_1$$
Hausdorff Young estimate ($1<p<2$) 
$$\|   \sqrt{\log |a(\infty)|} \|_{p'} \le C_p \|F\|_p$$
Plancherel identity
\begin{equation}\label{fourier-plancherel}
\|    \sqrt{\log |a(\infty)|} \|_{2}  = \frac \pi 2 \|F\|_2
\end{equation}
Maximal Riemann Lebesgue estimate
$$\|    \sup_x \sqrt{\log |a(x)|} \|_{L^\infty(k)} \le C \|F\|_1$$
Maximal Hausdorff Young estimate ($1<p<2$)
$$\|   \sup_x \sqrt{\log |a(x)|} \|_{L^{p'}(k)} \le C_p \|F\|_p$$
\end{theorem}
The Riemann Lebesgue and maximal Riemann Lebesgue estimates
follow easily from Gronwall's inequality, i.e., from applying
operator norms to (\ref{equivalent-ode}) and integrating the
inequality
$$\|G\|'/\|G\|\le \|W\|\ \ \ .$$
Then one uses $\sqrt{\log|a|}\sim \log\|G\|$ for small values
of $a$ and $\log|a|\sim \log\|G\|$ for large values
of $a$. We remark that in the $L^1$ theory one may view
$\log\|G\|$ as the more natural measure of the size of $G$
than $\sqrt{\log|a|}$.

The Hausdorff Young and
maximal Hausdorff Young inequalities follow by the work
of Christ and Kiselev \cite{ck0},\cite{ck1}. The Plancherel 
identity is a well known
scattering identity. Variants of it appear in \cite{buslaevfaddeev}
and \cite{verblunsky1},\cite{verblunskyII}. 
For the convenience of the reader
and to contrast it to our results in the $d$-adic model we will
sketch a proof in the appendix. Interestingly, while Plancherel gives
the $L^2$ endpoint of the Hausdorff Young inequality, we do not
known whether the constant $C_p$ in the Hausdorff Young
inequality can be chosen uniformly as $p$ tends to $2$.

The maximal version of Plancherel, which amounts to a scattering
variant of Carleson's theorem, is not known. We state
it as a conjecture

\begin{conjecture}
Carleson-Hunt estimate
$$\|   \sup_x \sqrt{\log |a(x)|} \|_{L^{2}(k)} \le C \|F\|_2\ \ .$$
\end{conjecture}
A more modest conjecture is 
\begin{conjecture}\label{weak-carleson-conjecture}
Weak type Carleson estimate
$$|\{k: \sup_x \sqrt{\log |a(x)|}>\lambda\}| \le 
C \lambda^{-2} \|F\|_2^2\ \ .$$
\end{conjecture}
Even more modestly one could conjecture that the function
$\sup_x \sqrt{\log |a(x)|}$ is finite almost everywhere for 
$F$ in $L^2(\R)$. To make $G$ well defined for this last conjecture
which is formulated in terms of the scattering transform for arbitrary
$F\in L^2(\R)$ one may replace the initial condition $G(-\infty)=\id$
by $G(0)=\id$. 

The main purpose of the current article is to give some supporting
evidence at least for Conjecture \ref{weak-carleson-conjecture}
by proving a variant of it in a $d$-adic model. 
The $d$-adic model is obtained by
replacing the exponential functions in (\ref{nonlinfour}),
which are the characters on $\R$, 
by characters of an infinite product of copies of $\Z(d)$ for some 
integer $d>1$. We call these groups Cantor groups.

From now on, $x$ and $k$ will denote non-negative real numbers.
For almost all such numbers, we have unique expansions with base $d$:
$$
k=\sum_{n\in \Z} k_n d^n
, \ \ \ 
x=\sum_{n\in \Z} x_n d^n
$$
where
$k_n$ and $x_n$ take values in $0,1,\dots, d-1$ and they are zero
for sufficiently large positive index $n$. Indeed, we shall
make these expansions unique for all $x$ and $k$ by requiring
each of them to have only finitely many non-zero entries
whenever possible.

Then we define a character function $w$ on $\R^+_0\times \R^+_0$ as
\begin{equation}\label{definewalsh}
w(k,x)=\gamma^{\sum_{n\in \Z} k_n
x_{-1-n}}\
\ ,
\end{equation} 
where $\gamma\in \C$ is some fixed primitive $d$-th root of unity.
The exact choice of $\gamma$ is not important.
Observe that the formally infinite sum in the exponent in (\ref{definewalsh})
has only finitely many non-zero summands.

Let $F\in L^2(\R^+)$.
For every parameter $k$ we consider the following initial value
problem:

\begin{equation}\label{wode}
\partial_x G(k,x)= W(k,x) G(k,x)
\end{equation}
$$G(k,0)={\rm id}$$
where
$$W(k,x)=\left(\begin{array}{cc} 0  & {F(x) w(k,x)} 
\\ \overline{F(x) w(k,x)} & 0 \end{array}\right)$$

By standard ODE theory this initial value problem has a unique
absolutely continuous solution satisfying the ODE almost everywhere.

We denote again by $a(k,x)$ the upper left entry of
$G(k,x)$.
The main theorem of this article is the following
\begin{theorem}\label{main}
Let $d>1$ and $F\in L^2(\R^+)$, and let $G$ be defined by (\ref{wode}).
Then for almost all $k\in \R^+$ the limit
\begin{equation}\label{limit}
G(k,\infty)=\lim_{x\to \infty} G(k,x)
\end{equation}
exists and satisfies the estimate
\begin{equation}\label{nlwp}
\int_0^\infty \log|a(k,\infty)|\, dk\le
C \int_0^\infty |F(x)|^2\, dx
\ \ \ .
\end{equation}
Moreover, 
\begin{equation}\label{nlwc}
|\{k: \sup_x |a(k,x)|>\lambda\}|\le C\lambda^{-1} 
{\|F\|_2^2}
\end{equation}
for all $\lambda>0$. Here as well as in (\ref{nlwc}) 
the constant $C$ may grow polynomially in $d$ but is
independent of $F$ and $\lambda$.
\end{theorem}

If $F$ is real valued, then a special situation occurs in 
Theorem \ref{main} for $d=2$: the matrices $W(k,x)$ then are real valued 
and commute for different values of $x$. By simultaneously diagonalizing 
all these matrices
one can decouple the two equations and obtain an ODE of the form
$\tilde{G}'=V\tilde{G}$ with
$$V(k,x)=\left(\begin{array}{cc}  F(x) w(k,x) & 0 
\\ 0 & - F(x) w(k,x) \end{array}\right)\ \ .$$
The solution at $+\infty$ of the corresponding initial value problem
is
$$V(k,x)=\left(\begin{array}{cc}  \exp(\hat{F}(x))& 0 
\\ 0 & \exp(- \hat{F}(x)) \end{array}\right)$$
where $\hat{F}$ denotes the Walsh-Fourier transform
(the Fourier transform with respect to the Cantor group with $d=2$).
In this special case Theorem \ref{main} follows simply from
the known Plancherel identity and Carleson's theorem
for the Walsh Fourier transform \cite{billard}. This example shows nicely
the connection of scattering transforms and the Fourier transform.

We will outline the proof of (\ref{nlwp}) in Section 
\ref{plancherel-section}. The proof is based on certain
swapping inequalities, which are discussed in detail
in Section \ref{swapping-section}. The proof of these
inequalities seems to be a genuinely new ingredient in the $d$-adic 
model as compared to the theory of the linear Fourier transform. 
In Section \ref{carleson-section} we prove (\ref{nlwc}),
which then easily implies (\ref{limit}).

Initially the authors had attempted to use multilinear
expansions of the solutions to  (\ref{wode}) to prove
Conjecture \ref{weak-carleson-conjecture} in the way 
Christ and Kiselev prove their results for $p<2$. 
However, as was observed in \cite{coex}, the terms in
this expansion do not satisfy reasonable bounds for $F\in L^2(\R)$.
Since the purpose of the current article is to compare
the $d$-adic to the continuous case, we prove a result
(Proposition \ref{unbounded-3})
in Section \ref{multilinear-section} which shows that the 
multilinear terms in the $d$-adic setting are equally badly
behaved. This is the second new result of this article.

In the appendix (Section \ref{appendix-section}) we sketch a proof
of the Plancherel identity (\ref{fourier-plancherel})
in Theorem \ref{classical}.
We only know a proof of this identity using complex contour integration.
This proof seems to not have the same flexibility as the proof
in the $d$-adic case which decomposes the scattering transform
into its elementary pieces. This in a sense is the main reason why at this
point we are unable to 
prove Carleson's theorem for the continuous scattering transform.

The first author was supported by NSF grant DMS 0100796.
The second author is a Clay Prize Fellow and is supported by a grant from
the Packard Foundations. The third author was supported
by a Sloan Fellowship and by NSF grants DMS 9985572 and DMS 9970469.

\section{Proof of the Plancherel inequality (\ref{nlwp})}\label{plancherel-section}

First we consider
the case of compactly supported $F$. Thus, for fixed $k$,
$G(k,x)$ becomes constant for large $x$ and the existence of
the limit $G(k,\infty)$ is not in question.

Recall that $SU(1,1)$ is the Lie group of all complex 
$2\times 2$ matrices of the form 
\begin{equation}\label{su11}
\left(\begin{array}{cc} a & b \\ \overline{b} &
\overline{a}\end{array}\right)
\end{equation} 
with determinant $|a|^2-|b|^2=1$. This
group is isomorphic to $SL_2(\R)$. Observe that $W(k,x)$ is an element of
the Lie algebra of $SU(1,1)$, and thus the solution to the
initial value problem (\ref{wode}), which is well known to exist as an 
absolutely 
continuous function, takes values in $SU(1,1)$.
Of course one can verify directly by an elementary calculation
that the solution to
(\ref{wode}) has the form (\ref{su11}) and
determinant $1$ for all $x$, which is all we need from this brief 
discussion of Lie groups.

The following is an easy observation
about breaking the ODE (\ref{wode}) into pieces along the $x$ variable.
For any interval $\omega\subset \R^+_0$ define the localized system

\begin{equation}\label{localized}
\partial_x G_\omega(k,x)=W_\omega(k,x) G_\omega(k,x)
\end{equation}
$$G_\omega(k,0)=\id\ \ \ ,$$
where
$$W_\omega(k,x)=\left(\begin{array}{cc} 0 & F(x) 1_{\omega}(x)w(k,x) 
\\ \overline{F(x) 1_{\omega}(x) w(k,x)} & 0 \end{array}\right)\ \ .$$

\begin{lemma}\label{partition}
Let $\omega_1, \omega_2, \dots, \omega_n$ be adjacent intervals in ascending
order, and let the union of these intervals be the interval $\omega$.
Then we have  for all $k\ge 0$:
$$G_\omega(k,\infty)= \prod_{j=n}^1 G_{\omega_j}(k,\infty)=
G_{\omega_n}(k,\infty)\dots G_{\omega_2}(k,\infty)G_{\omega_1}(k,\infty)\ \ .$$
\end{lemma}

Proof: By induction the lemma follows from the special case for
two adjacent intervals $\omega_1$ and $\omega_2$.
Fix $k$. It is easy to check that the absolutely continuous function
\begin{equation}\label{matrixproduct}
G_{\omega_2}(k,x)G_{\omega_1}(k,x) 
\end{equation}
satisfies the differential equation for $G_{\omega}(k,x)$ almost everywhere.
This follows easily from letting $x_0$ be the point separating 
$\omega_1$ and $\omega_2$ and considering $x<x_0$ and $x>x_0$ separately.
Since (\ref{matrixproduct}) also satisfies the correct initial
condition, this proves the lemma.

\endproof

Next, we claim that if $\omega$ is a $d$-adic interval, that means
an interval of the form
$$[d^\kappa n,d^\kappa (n+1)) $$
with integers $\kappa$ and $n\ge  0$,
then $G_\omega(k,x)$ does not change much as $k$ varies
inside a $d$-adic interval of reciprocal length $d^{-\kappa}$.

To make this claim precise, we define a tile to be a 
rectangle 
$p=I\times \omega$  of the form
$$[d^\kappa n,d^\kappa (n+1))\times [d^{-\kappa} l , d^{-\kappa} (l+1))$$
with integers $\kappa,n,l$ such that $l,n\ge 0$.

\begin{lemma}\label{constance}
Let $I\times \omega$ be a tile. Let $k_0$ be the left endpoint of
$I$ and let $k$ be any point in $I$.
Then there is an integer $j=j(k)$ independent of $x$ such that if
$$G_\omega(k_0,x)=\left(\begin{array}{cc} a & b \\ \overline{b} & \overline{a}
\end{array}\right)\ \ ,$$
then 
$$G_\omega(k,x)=\left(\begin{array}{cc} a & { \gamma^{j} b}\\ 
\overline{\gamma^{j} b} & 
\overline{a}
\end{array}\right)\ \ .$$
In particular, the first entry $a$ of $G_\omega(k,x)$
is independent of $k$ as long as $k\in I$. 
\end{lemma}

Proof:
Assume the length of $I$ is $d^\kappa$.
Since $G_\omega$ is constant outside $\omega$, it suffices to
show the claim for $x\in \omega$.
We split $w(k,x)$ into two factors as follows:
$$w(k,x)=
(\gamma^{\sum_{\nu < \kappa} k_\nu
x_{-1-\nu}})
(\gamma^{\sum_{\nu \ge \kappa} k_\nu
x_{-1-\nu}})\ \ \ .$$
Observe that if $x$ varies in $\omega$, then the first
factor in this splitting does not change. Likewise,
the second factor is constant for $k\in I$. Thus there
is a $j$ depending on $k\in I$ such that
$$w(k,x)1_{\omega}(x)= \gamma^j w(k_0,x)1_{\omega}(x)\ \ .$$
Now let $\Gamma$ be the constant matrix 
$$\Gamma =
\left(\begin{array}{cc} \gamma^{j} & 0\\ 
0 & 1  
\end{array}\right)\ \ \ .$$
Then
$$W_{\omega}(k,x)= \Gamma\, W_{\omega}(k_0,x)\, \Gamma^{-1}\ \ .$$
By conjugating the initial value problem (\ref{localized}) by $\Gamma$
we observe that
$$G_\omega(k,x)= \Gamma\, G_{\omega}(k_0,x)\, \Gamma^{-1}\ \ .$$
This proves the lemma.
\endproof
Motivated by this lemma we shall define for a tile $p=I\times \omega$:
$$G_p=G_\omega(k_0,\infty)$$
where $k_0$ is the left endpoint of $I$.
Next, we shall investigate the relation of the matrices $G_p$
for nearby tiles $p$. Here we mean by nearby tiles that the tiles
are contained in a given
$d$-adic rectangle of area $d$.

Define a multitile to be a rectangle $P=I\times \omega$
of the form 
$$[d^\kappa n,d^\kappa (n+1))\times [d^{1-\kappa} l , d^{1-\kappa} (l+1))$$
with integers $\kappa,n,l$ and $l,n\ge 0$.
There are $d$ tiles $p_j$, $j=0,\dots, d-1$
contained in $P$ of the form $I\times \omega_j$. We shall 
always assume the $\omega_j$ are ordered in ascending order. 
We call these tiles the horizontal subtiles of $P$.
Moreover, there are $d$ tiles $q_j$, $j=0,\dots, d-1$ contained
in $P$ of the form $I_j\times \omega$. We shall again
assume the $I_j$ are ordered in ascending order, and we call these
tiles the vertical subtiles of $P$.

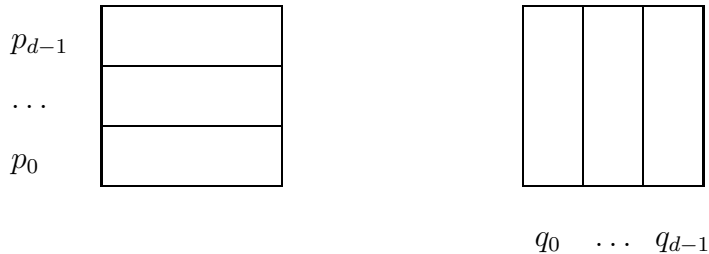
\setlength{\unitlength}{0.8mm}
\begin{picture}(160,55)

\put(40,15){\line(1,0){30}}
\put(40,25){\line(1,0){30}}
\put(40,35){\line(1,0){30}}
\put(40,45){\line(1,0){30}}

\put(40,15){\line(0,1){30}}
\put(70,15){\line(0,1){30}}

\put(25,18){$p_0$}
\put(25,28){$\dots$}
\put(25,38){$p_{d-1}$}

\put(110,15){\line(0,1){30}}
\put(120,15){\line(0,1){30}}
\put(130,15){\line(0,1){30}}
\put(140,15){\line(0,1){30}}

\put(110,15){\line(1,0){30}}
\put(110,45){\line(1,0){30}}

\put(112,5){$q_0$}
\put(122,5){$\dots$}
\put(132,5){$q_{d-1}$}

\end{picture}

\begin{lemma}\label{swappingid}

Let $P=I\times \omega$ be a multitile and assume
its horizontal tiles are $p_0,\dots,p_{d-1}$ and its
vertical tiles are $q_0,\dots,q_{d-1}$.
If
$$G_{p_j}=
\left(\begin{array}{cc} a_j & b_j\\ 
\overline{b_j} & 
\overline{a_j}
\end{array}\right) $$
for $j=0,\dots,d-1$, then
$$G_{q_m}= \prod_{j=d-1}^{0} 
\left(\begin{array}{cc} a_j & { \gamma^{mj} b_j}\\ 
\overline{\gamma^{mj} b_j} & 
\overline{a_j}
\end{array}\right)\ \ .$$
Here the product is to be read in descending order
$$\left(\begin{array}{cc} a_{d-1} & { \gamma^{m(d-1)} b_{d-1}}\\ 
\overline{\gamma^{m(d-1)} b_{d-1}} & 
\overline{a_{d-1}}
\end{array}\right)
\dots
\left(\begin{array}{cc} a_1 & { \gamma^{m} b_1}\\ 
\overline{\gamma^{m} b_1} & 
\overline{a_1}
\end{array}\right)
\left(\begin{array}{cc} a_0 & { b_0}\\ 
 \overline{b_0} & 
\overline{a_0}
\end{array}\right)\ \ .$$
\end{lemma}

Proof:
Let $k$ denote the left endpoint of $I$ and let $d^\kappa$
be the length of $I$. Let $p_j= I\times \omega_j$.
By Lemma \ref{partition} it suffices to prove for $0\le m\le d-1$:
\begin{equation}\label{j-th-matrix}G_{\omega_j}(k+d^{\kappa-1}m,\infty)= 
\left(\begin{array}{cc} a_{j} & { \gamma^{mj} b_{j}}\\ 
\overline{\gamma^{mj} b_{j}} & 
\overline{a_{j}}
\end{array}\right)\ \ .
\end{equation}
However, we have for $x\in \omega_j$:
$$w(k+d^{\kappa-1}m, x)=
(\gamma^{\sum_{\nu < \kappa-1} k_\nu
x_{-1-\nu}})
\gamma^{m j}
(\gamma^{\sum_{\nu \ge \kappa} k_\nu
x_{-1-\nu}})$$
$$= \gamma^{mj} w(k,x)
$$
Now (\ref{j-th-matrix}) follows by the considerations
in the proof of Lemma \ref{constance}.
\endproof

In the next section we will obtain a
function $\beta:SU(1,1)\to \R^+_0$ such that
$$\frac 1C \log|a|\le \beta(G)\le C \log|a|$$ 
for some constant $C$ depending polynomially on $d$
and,
with the notation of Lemma \ref{swappingid},
\begin{equation}\label{sw}
\sum_{m=0}^{d-1}\beta(G_{p_m})\le d \sum_{j=0}^{d-1} \beta(G_{q_j})\ \ .
\end{equation}
We will refer to $\beta$ as the {\it swapping function} and 
(\ref{sw})
as the {\it swapping property}. Assume for now this swapping function 
has been constructed. The rest of this section is to prove
(\ref{nlwp}) using this function.

Let $K$ be a large integer and consider the rectangle
$R=[0,d^K)\times [0,d^K)$. Let $\p_\kappa$ denote
the set of all tiles $I\times \omega\subset R$
with $|I|=d^\kappa$. One can partition the tiles in $\p_\kappa$
into $d$-tuples such that each $d$-tuple consists of the horizontal
tiles of a multitile.
Applying (\ref{sw}) on each tuple we obtain
$$\sum_{p\in \p_\kappa} \beta(G_p)\le d \sum_{p\in \p_{\kappa+1}} \beta(G_p)\ \ \ .$$
By iterating this we obtain
$$d^{-K} \sum_{p\in \p_{-K}} \beta(G_p)\le d^{K} \sum_{p\in \p_{K}} \beta(G_p)\ \ \ .$$
Hence
\begin{equation}\label{loga}
d^{-K} \sum_{p\in \p_{-K}} \log |a_p|\le C d^{K} \sum_{p\in \p_{K}} \log|a_p|
\end{equation}
where $a_p$ denotes the upper left entry of $G_p$. Observe that 
since we have no control over $K$ it is very important that there
are no further constants on the right hand side of (\ref{swappingid})
other than the constant $d$ which is the natural scaling constant
(as we will see momentarily).

We may assume that the support of $F$ is contained in $[0,d^K)$.
Then for a tile $p=I\times [0,d^K)$ we have that
$a_p$ is equal to $a(k,\infty)$ where $k$ is the left 
endpoint of $I$, or, by Lemma \ref{constance} where $k$ is
any point in $I$. Thus the left hand side of (\ref{loga})
is equal to
$$\int_{0}^{d^K} \log |a(k,\infty)| \, dk \ \ \ .$$

Thus it remains to show that the right hand side of (\ref{loga})
is less than 
$$C \|F\|_2^2 $$
for arbitrarily large $K$ and constant $C$ independent of $K$.
Observe that (\ref{localized}) implies
$$\frac {\partial}{\partial x}\|G_\omega\|_{op} \le \|W_\omega\|_{op}
\|G_\omega\|_{op}$$
which implies together with the intitial condition for $G$ at $0$:
$$ \|G_\omega(k,\infty)\|_{op}\le \exp(\int_0^\infty \|W_\omega\|_{op})\ \ \ .$$
This is Gronwall's inequality and - as has been mentioned before -
implies the $L^1$  estimates claimed in Theorem \ref{classical}.
Continuing the present considerations we obtain 
$$ \log \|G_\omega(k,\infty)\|_{op}\le \|F\|_{L^1(\omega)}\ \ \ .$$
We claim that the operator norm of a matrix of the type (\ref{su11})
is equal to $|a|+|b|$.
This is clear in the case that
$a$ and $b$ are real, in which it is easy to calculate the eigenvalues of the
symmetric matrix. The general case can be obtained by multiplying 
the matrix from both sides by
unitary diagonal matrices to reduce to the real case.

By H\"older's inequality we thus have
$$\log(|a_\omega|+|b_\omega|)\le \|F\|_2 |\omega|^{1/2}\ \ \ .$$
By choosing $\omega$ small enough ($K$ large enough), the right hand side
can be made small. Thus we can assume $|a|$ is close
to $1$ and $b$ is close to $0$. Then we obtain
$$ \log (|a_\omega|+|b_\omega|)\ge 
\frac 12 |b_\omega|\ge \frac 14 \sqrt{\log|a_\omega|}\ \ .$$
Hence
$$  d^{K} \sum_{p\in \p_{K}} \log|a_p| \le 
4 \sum_{p\in \p_{K}} \|F_\omega\|_2^2 = 4 \|F\|_2^2\ \ .$$
This gives the desired bound on the right hand side
of (\ref{loga}) and completes the proof of inequality (\ref{nlwp})
in the case of compactly supported $F$.

If $F$ is not compactly supported then we will show later
that the limit (\ref{limit}) exists almost everywhere. 
Assuming this for now, then (\ref{nlwp}) follows by 
Fatou's lemma.

\section{The swapping function}\label{swapping-section}

In this section we will find for each $d>1$ a function
$$\beta:SU(1,1)\to \R^+_0$$ such that $\beta(G)$ is comparable
to $\log|a|$ and we have the inequality (\ref{sw}).
We call this inequality a swapping inequality and $\beta$
a swapping function because (\ref{sw}) swaps the vertical
tiles to the horizontal tiles in a given multitile.

The case $d=2$ is particularly easy and we will do
it first. The function $\beta$ simply can be chosen
to be the logarithm of the Hilbert Schmidt norm of $G$.
Here we define the Hilbert Schmidt norm of
a matrix $G$ of the form (\ref{su11}) to be
$$ \|G\|_{HS}=\sqrt{|a|^2+|b|^2}\ \ \ .$$
Observe that for $a$ near $1$ we have
\begin{equation}\label{small-a}
\log\|G\|_{HS}\sim  \log (|a|^2+|b|^2) \sim |b|^2 \sim \log|a|
\end{equation}
and for large $a$ we have
\begin{equation}\label{large-a}
\log\|G\|_{HS}\sim  \log (2 |a|^2 ) \sim \log|a|\ \ \ .
\end{equation} 
Thus $\beta(G)$ is comparable to $\log|a|$.
We write
\begin{equation}\label{minusstar}
A^{-*}:=(A^{-1})^*
\end{equation}
The following lemma then says that the swapping inequality is true.
\begin{lemma}\label{swapping}
If $A, B\in SU(1,1)$ then
$$\log \|AB\|_{HS} +\log \|A B^{-*}\|_{HS} 
\le 2\log\|A\|_{HS}+2\log\|B\|_{HS}$$
\end{lemma}

Write
$$ A= \left(\begin{array}{cc} a & b \\ \overline{b} & \overline{a}
\end{array}\right)$$
and 
$$ B= \left(\begin{array}{cc} c & d \\ \overline{d} & \overline{c}
\end{array}\right),\ \ 
 B^{-1}= \left(\begin{array}{cc} \overline{c} & - d \\ - \overline{d} & {c}
\end{array}\right),\ \ 
B^{-*}= \left(\begin{array}{cc} c & - d \\ - \overline{d} & \overline{c}
\end{array}\right)\ \ \ .$$
Then we have
$$AB = \left(\begin{array}{cc} ac + b\overline{d}  
& b \overline{c}+a\overline{d}
\\
\overline{b}c+\overline{a}d &
\overline{a}\overline{c}+\overline{b}d
\end{array}\right),\ \ AB^{-*} = \left(\begin{array}{cc} ac - b\overline{d} 
& b \overline{c} - a\overline{d}
\\
\overline{b}c - \overline{a}d &
\overline{a}\overline{c} - \overline{b}d
\end{array}\right)\ \ .$$
This gives
$$\|AB\|_{HS}^2 + \|AB^{-*}\|_{HS}^2= 2(|ac|^2+|bd|^2+|bc|^2+|ad|^2)
= 2\|A\|_{HS}^2\|B\|_{HS}^2\ \ .$$
Using the arithmetic mean-geometric mean inequality we obtain
$$\|AB\|_{HS}\|AB^{-*}\|_{HS}\le \|A\|_{HS}^2\|B\|_{HS}^2\ \ \ .$$
Taking logarithms proves the lemma.
\endproof
We remark that the function $\beta(G):=\log|a|$ does not satisfy
the swapping inequality in general. This can be seen from choosing
$a$,$b$,$c$ positive and $d$ purely imaginary in the above example.

Now consider $d\ge 3$. In this case one has to choose a more
complicated swapping function. Indeed, in an appendix
to this section we will sketch an argument that the logarithm of the
Hilbert Schmidt norm does not satisfy the required 
swapping inequality.

Choose an $\epsilon$ sufficiently small. For the purpose
of keeping track of polynomial growth in the parameter
$d$ we remark that the choice 
$10^{-3}d^{-1}$ will be sufficient.

Let $r$ be the smallest positive number such that
$$r^2- r^3= \epsilon^{10}+\epsilon^{20} {\rm arcsinh}(r)\ \ .$$
Then $r$ is of the order $\epsilon^{5}$.
We consider the following swapping function defined on $\C$
$$\beta(z):= |z|^2 -  |z|^3$$
if $|z|\le r$ and 
$$\beta(z):=\epsilon^{10}+\epsilon^{20} {\rm arcsinh} (|z|)$$
if $|z|\ge r$.

Finding this function was inspired by the discussion of
Bellman functions in \cite{nt}, whence the letter $\beta$ for
this function. For a matrix $G$ of type (\ref{su11}) we
will let $\beta(G)=\beta(b)$. By a discussion as in (\ref{small-a})
and (\ref{large-a}) it is clear that $\beta(G)$ is comparable
to $\log|a|$ with constants growing polynomially in $\epsilon^{-1}$ and thus
growing polynomially in $d$.

Clearly there is not a unique way to choose $\beta$. Our choice
reflects in a very explicit way the two different types of 
behaviour for small and
for large $|z|$ which will be apparent from the discussion
below. Moreover, for each of the two regions our choice shows
explicitly the leading order term ($|z|^2$ and 
$C^{-1}{\rm arcsinh(|z|)}$) and a smaller order correction 
term which is used to estimate the nonlinear effects.
The third order correction term for small $z$ could be replaced
by any other power $|z|^p$ with $2<p<4$.

Given $d$ pairs $(a_i,b_i)$ of complex numbers
with $|a_i|^2=1+|b_i|^2$ and a $d$-th root of unity $\gamma$
(in this section $\gamma$ shall not be a fixed primitive $d$-th root
of unity but an arbitrary $d$-th root of unity) we define

$$
\left(\begin{array}{cc} A_\gamma & {B_\gamma}\\ \overline{B_\gamma} &
\overline{A_\gamma}
 \end{array}\right)
=\prod_{i=1}^d \left(\begin{array}{cc} a_i & {\gamma^{i}b_i}
\\ \overline{\gamma^{i}b_i} & \overline{a_i}
 \end{array}\right)\ \ .
$$
The factors in this product do not commute, hence we
emphasize that the product is understood in ascending order:
\begin{equation}\label{matrix-product} 
\left(\begin{array}{cc} a_1 &
{\gamma^{1}b_1}
\\ \overline{\gamma^{1}b_1} & \overline{a_1}
 \end{array}\right)
\left(\begin{array}{cc} a_2 & {\gamma^{2}b_2}
\\ \overline{\gamma^{2}b_2} & \overline{a_2}
 \end{array}\right)
\dots
\left(\begin{array}{cc} a_d & {\gamma^{d}b_d}
\\ \overline{\gamma^{d}b_d} & \overline{a_d}
 \end{array}\right)\ \ .
\end{equation}
Also observe that in the last factor we have $\gamma^d=1$.
 
\begin{lemma}\label{swapping-inequality-lemma}
Under the above hypotheses, we have
$$\sum_{\gamma} \beta(B_\gamma)\le d \sum_{i=1}^d \beta(b_i)$$
where the sum on the left hand side runs over all $d$-th
roots of unity.
\end{lemma}

This lemma clearly implies the desired swapping inequality (\ref{sw}).

Proof:
We shall first consider the case when $|b_i|\le r$
for all $i$. 

Observe that we can write $B_\gamma$ as a polynomial in the (for this matter
viewed as independent) variables $a_i$, $b_i$, $\overline{a_i}$
and $\overline{b_i}$ for $i=1,\dots,d$. 

We claim that this polynomial is odd in the vector  
$${\bf b}=(b_1,\overline{b_1},\dots,b_d,\overline{b_d})\ \ \ .$$ 
This follows from the observation that the operation
$G\to G^{-*}$ (see (\ref{minusstar}) and the lines thereafter) commutes with matrix products. 
Thus replacing the vector $\bf b$
by its negative replaces $B_\gamma$ by its negative. Thus $B_\gamma$
has to be an odd polynomial in $\bf b$. 

Writing down the matrix product explicitly, we observe that
the polynomial $B_\gamma$ is a sum of monomials of degree $d$,
where each such monomial has exactly one entry from each of the
matrices $G_1,\dots,G_{d}$ as factor. Any choice of one entry
from each matrix can appear in a monomial, provided the following 
row- and column conditions are satisfied: an entry from the $j$-th 
row of $G_i$ can appear only if an entry from the $j$-th column of 
$G_{i-1}$ appears, the entry form $G_{d}$ has to be from
the second column, and the entry from $G_{1}$ has to be from the
first row.

This together with oddness in $\bf b$ gives the crude estimate
$$|B_\gamma|\le \sum_i |b_i| \prod_{j\neq i}(|a_i|+|b_i|)
\le d r(1+2r)^{d-1}\le 2 d r\ \ \ .$$
The right hand side is of the order $\epsilon^4$ and thus
we are well in the range such that we have an
estimate
\begin{equation}\label{quadratic-range}
\beta(B_\gamma)\le |B_\gamma|^2 -
|B_\gamma|^3\ \ \ .
\end{equation}

Studying now the polynomial of $B_\gamma$ more carefully, we
extract those terms  which are linear in $\bf b$. They
are easily seen to be
\begin{equation}\label{linear-terms}
\sum_{i=1}^{d} \gamma^{i}
\left(\prod_{j<i}{a_j}\right) b_i
\left(\prod_{j>i}\overline{a_j}\right)\ \ .
\end{equation}
Moreover, since there are no terms quadratic in $\bf b$,
we obtain the estimate
\begin{equation}\label{linear-taylor}
\left|B_\gamma -\sum_{i=1}^{d} \gamma^{i}
\left(\prod_{j<i}{a_j}\right) b_i
\left(\prod_{j>i}\overline{a_j}\right)\right|
\end{equation}
$$\le \sum_{i< j < k}|b_i||b_j||b_k|
\prod_{l\neq i,j,k}|a_l|+|b_l| \le 
2 d^3 |b_m||b_{m'}|^2\ \ .
$$
Here $m$ denotes the index such that $|b_m|$ is maximal
among all $|b_i|$ and $m'$ denotes the index such that $|b_{m'}|$
is maximal among all $b_i$ with $i\neq m$.

Now we consider the polynomial for $|B_\gamma|^2=B_\gamma\overline{B_\gamma}$ 
and sum
over all $\gamma$.  
Observe that upon the  summation in $\gamma$, all terms of the 
polynomial of $|B_\gamma|^2$
which have a non-trivial power of $\gamma$ in the coefficient get canceled.

We are again interested in the lowest order terms in $\bf b$, which
are the quadratic terms.
Such terms appear 
when the $i$-th summand of (\ref{linear-terms}) is multiplied
by the complex conjugate of the $j$-th summand in (\ref{linear-terms}).
The power of $\gamma$ in such a term is trivial only if $i=j$. 
Thus the (in $\bf b$) quadratic terms of $\sum_\gamma |B_\gamma|^2$
are precisely
\begin{equation}\label{quadratic terms}    
d \sum_{i=1}^{d} |b_i|^2
\left(\prod_{j\neq i}|a_j|^2\right)\ \ \ .
\end{equation}
Moreover, from the previous discussion 
we can easily see the estimate
$$\left|\sum_\gamma |B_\gamma|^2 -d \sum_{i=1}^{d} |b_i|^2
\left(\prod_{j\neq i}|a_j|^2\right) \right|
\le 4 d^4 r |b_m||b_{m'}|^2\ \ \ .
$$
Now we use the fact that $|a_i|^2=1+|b_i|^2$
to obtain
\begin{equation}\label{quadr-term}
\left|\sum_\gamma |B_\gamma|^2 -d \sum_{i=1}^{d} |b_i|^2 \right|
\le 6 d^4 r |b_m||b_{m'}|^2\ \ \ .
\end{equation}
Next, we observe 
$$ 
\frac 1d \sum_{\gamma} |B_\gamma|^3\ge 
\left(\frac 1d \sum_{\gamma} |B_\gamma|  ^2\right)^{3/2}
\ge\left( \sum_{i=1}^3|b_i|^2\right)^{3/2}-\frac 32 
6 d^3 r |b_m||b_{m'}|^2\ \ \ .
$$ 
The last estimate followed from (\ref{quadr-term}) and a trivial
estimate on the slope
of the function $x\to x^{3/2}$ in the interval $[0,1]$. Now the right hand 
side of the last display is equal to
$$|b_m|^3 \left( 1+ \sum_{i\neq m}(|b_i|/|b_m|)^2\right)^{3/2}- 
9 d^3 r |b_m||b_{m'}|^2
$$
$$\ge |b_m|^3 \left( 1+ \frac 32 \sum_{i\neq m}
(|b_i|/|b_m|)^2\right)-9 d^3 r |b_m||b_{m'}|^2
$$
$$\ge \sum_i |b_i|^3 +  \frac 12 |b_m||b_{m'}|^2 -
9 d^3 r |b_m||b_{m'}|^2
\ge \sum_i |b_i|^3 +  \frac 14 |b_m||b_{m'}|^2\ \ .
 $$

Together with (\ref{quadr-term}) we obtain
$$\sum_{\gamma} |B_\gamma|^2 -  |B_\gamma|^3
\le d \sum_{i} |b_i|^2 -  |b_i|^3\ \ .$$
Together with   (\ref{quadratic-range})
this proves the Lemma \ref{swapping-inequality-lemma}
in the present case $|b_i|\le r$ for all $i$.

Now we consider the case at the other extreme that for at least two
indices $i$ we have  $|b_i|\ge r$. Denote by $I$ the set of indices $i$ for which
$|b_i|\ge r$.

We observe that the quantity $\arcsinh|b|$ has the meaning of the logarithm
of the operator norm of the matrix
$$\left(\begin{array}{cc} a & {b}
\\ \overline{b} & \overline{a}
 \end{array}\right)$$
By elementary calculus using $|a|^2= 1+ |b|^2$ this is equivalent to the statement that
the operator norm of this matrix is $|a|+|b|$. This however has been observed
in Section \ref{plancherel-section}.

In particular we obtain for each $\gamma$:
$$\arcsinh|B_\gamma|\le \sum_i \arcsinh |b_i|\ \ \ .$$
Since by elementary calculus we always have
$$\beta(B_\gamma)\le \epsilon^{10} + \epsilon^{20} \arcsinh|B_\gamma|\ \ \ .$$
we obtain
$$\beta(B_\gamma)\le \epsilon^{10} + \epsilon^{20} \sum_i \arcsinh|b_i|$$
$$\le \epsilon^{10}(1-|I|) + \sum_{i\in I} \beta(b_i) + \epsilon^{20} \sum_{i\notin I} \arcsinh|b_i|$$
$$\le - \epsilon^{10} + \sum_{i\in I} \beta(b_i) + \epsilon^{20} dr
\le \sum_{i\in I} \beta(b_i)\ \ .$$
Summing over $\gamma$ proves Lemma \ref{swapping-inequality-lemma}
in this case.

The same reasoning as in the previous case can be applied
if there is only one index $i$ such that $|b_i|\ge r$
but there is at least one other index $j$ such that $|b_j|>\epsilon r$.
The latter implies that
$$\epsilon ^{20} \arcsinh|b_j| \le \beta(b_j) - \epsilon ^{15}\ \ .$$
Namely, the left hand side is less than $\epsilon^{20}$, while 
$\beta(b_j)$ is at least $\epsilon^{14}$. Thus we have
$$\beta(B_\gamma) \le \beta(b_i) + \beta(b_j) -\epsilon^{15} + 
\epsilon^{20} \sum_{k\neq i,j}\arcsinh |b_k|$$
$$\le \beta(b_i) + \beta(b_j)\ \ .$$
This proves the Lemma \ref{swapping-inequality-lemma} in the given case.

It remains to prove the case when there is one index $i$
such that $|b_i|\ge r$ and for all other indices $j\neq i$
we have $|b_j|\le \epsilon r$. We extract from (\ref{linear-taylor}):
$$\gamma^{-i}B_\gamma = b_i + \sum_{j>i} \gamma^{j-i}a_i b_j + 
\sum_{j<i}\gamma^{j-i}\overline{a_i}b_j +E$$
with
$$|E|\le 4d^2(|a_i|+|b_i|)|b_{m'}|^2
 \le 8d^2(1+|b_i|)|b_{m'}|^2\ \ .$$
Here $m'$ is the index such that $|b_{m'}|$ is maximal among
all $|b_j|$ with $j\neq i$.
Observe that under the given assumptions the term $b_i$ is large
compared to the linear terms in $b_j$, $j\neq i$, which in
turn are large compared to $E$. Indeed, we observe the estimate
$$\left|\sum_{j<i} \gamma^{j-i}a_i b_j + 
\sum_{j>i}\gamma^{j-i}\overline{a_i}b_j +E \right|\le 
2d(1+|b_i|) |b_{m'}|\ \ .$$

Our goal is to make a Taylor expansion of the function
$f:z\to \arcsinh|z|$ near the point $b_i$. Let
$\lambda$ denote the linear form which is the
derivative of $f$ at $b_i$ and let $\rho$ denote
the quadratic form which is the second derivative of
$f$ at an appropriate point within distance
$2d(1+|b_i|) |b_{m'}|$ of $b_i$.

Then we obtain from Taylor's theorem
$$\sum_\gamma \arcsinh|B_\gamma| =\sum_\gamma \arcsinh|\gamma^{-i}B_\gamma|$$
$$= d \, \arcsinh|b_i|+ \sum_\gamma \lambda( \sum_{j< i} \gamma^{j-i}a_i b_j + 
\sum_{j>i}\gamma^{j-i}\overline{a_i}b_j +E)+F$$
with $$|F|\le \sum_\gamma |\rho||\sum_{j< i} \gamma^{j-i}a_i b_j + 
\sum_{j>i}\gamma^{j-i}\overline{a_i}b_j +E|^2\ \ .$$
Using that $$\sum_\gamma \gamma^{j-i}=0$$ for all $j\neq i$ we obtain
$$\sum_\gamma \arcsinh|B_\gamma| =d\,  \arcsinh|b_i| + \sum_\gamma \lambda(E) 
+ F\ \ .$$
From elementary calculus we obtain
$$ |\lambda|\le (1+|b_i|^2)^{-1/2}$$
$$ |\rho|\le 20 |b_i|^{-1}(1+|b_i|^2)^{-1/2}\ \ .$$

Hence
$$\sum_\gamma \arcsinh|B_\gamma| \le d\,  \arcsinh|b_i| +  50 r^{-1} d^3 
 |b_m|^2$$
and 
$$\sum_\gamma \beta(B_\gamma) \le d \beta(b_i) + \epsilon^{20} 50  r^{-1} d^3  
 |b_m|^2
\le d \sum_{j} 
 \beta(b_j)\ \ .$$
This proves the Lemma \ref{swapping-inequality-lemma} in the last case and thus completes the proof.

\endproof

We close this section by showing that for $d=3$ we cannot choose
$\beta(G)$ to be $\log\|G\|_{HS}$ if we want the swapping 
inequality (\ref{sw}) to hold. The rest of this section is irrelevant for the 
purpose of proving Theorem \ref{main}.

Let $\gamma$ denote a third root of unity.
Define
$$
\left(\begin{array}{cc} a & b\\ c & d
 \end{array}\right)_\gamma
=
\left(\begin{array}{cc} a & \gamma b\\ \gamma^2 c & d
 \end{array}\right)
$$

We aim to find matrices $A,B,C$ in
$SU(1,1)$ such that
\begin{equation}\label{counterexample}
\prod_\gamma \| A_{\gamma \gamma}  B_\gamma C \|_2
> \|A\|^3\|B\|^3\|C\|^3
\end{equation}
where the product on the left hand side goes
over all three third roots of unity.
We shall present such matrices $A$, $B$, $C$
with real entries. Thus we write
$$
A=\left(\begin{array}{cc} a & b\\ b & a
 \end{array}\right),\ \ 
B =
\left(\begin{array}{cc} c  & d\\ d & c
 \end{array}\right),\ \ 
C= \left(\begin{array}{cc} e & f\\ f & e
 \end{array}\right)\ \ .
$$ 
By homogeneity of (\ref{counterexample})
the requirement that $A,B,C$ are in $SU(1,1)$
can be relaxed to the requirement that they have
nonzero determinant and $|a|>|b|$, $|c|>|d|$, and
$|e|>|f|$. Indeed we will produce an example satisfying
the latter constraints and 
$$a^2+b^2=1,\ \ c^2+d^2=1,\ \ e^2+f^2=1\ \ \ .$$
In particular the right hand side of (\ref{counterexample}) 
is equal to $1$.

The matrix $A_{\gamma\gamma} B_\gamma C$ is equal to
$$
A=\left(\begin{array}{cc} ace + \gamma a{d}f+ \gamma bde+ \gamma^2
b{c} f & *
\\ 
\gamma bce + \gamma^2 b {d} f + \gamma^2 ade + a{c} f
 & *
 \end{array}\right)$$
where the unspecified entries in the second column are the same
as the diagonally opposite terms.

We calculate the Hilbert Schmidt norm squared of this matrix,
which is the sum of modulus sqared of the two indicated entries..
Observe that squaring the entries and multiplying out
gives pure squares  (the modulus square of a summand)
and mixed terms (product of two different summands with proper
complex conjugation). The pure squares simply add up to
$$(a^2+b^2)(c^2+d^2)(e^2+f^2)=1\ \ .$$
To calculate the mixed terms it helps to observe that we may
divide the second entry by $\gamma$, then the two entries
are alike but with $a$ and $b$ interchanged. 
Using $a^2+b^2=1$ and $\overline{\gamma}=\gamma^2$ we
obtain for the mixed terms
$$ + cd  ef  \gamma^2
 +   2 ab c {d} e^2 \gamma^2 
+ 2ab c^2 ef \gamma
+ 2ab {d}^2 ef
+ 2ab  c{d} f^2 \gamma^2
+ cd ef \gamma^2
$$
$$ + {c} {d} ef  \gamma 
+ 2 ab {c} d e^2 \gamma 
+ 2ab {c}^2 ef  \gamma^2
+ 2ab  {d}^2 ef
+ 2ab  {c}d f^2 \gamma
+ {cd} ef \gamma\ \ .
$$
Now we set $\alpha=ab$ and $\beta=ef$. 
Thus it will suffice to produce $\alpha,\beta\in
[0,{1/2}]$. Using $e^2+f^2=1$ we obtain for the square of
the Hilbert Schmidt norm
$$(1+ 4 \alpha\beta d^2)
+ ( 2\alpha\beta c^2+ 2\alpha {c} d + 
2\beta {c}{d})\gamma
+ ( 2\alpha\beta{c}^2+ 2\alpha c{d}  + 
2\beta {c}{d})\gamma^2\ \ .$$
Now we observe for any three numbers $K,L,M$ the formula
$$\prod_\gamma K+ L \gamma  + M \gamma^2 = K^3+L^3+M^3-3KLM\ \ .$$
Namely, expanding the left hand side, clearly
the coefficients in front of $K^3$, $K^2 L$, and $K^2 M$
are $1$, $0$, $0$; the latter two because the sum of
the third roots of unity is $0$. By multiplying each factor
on the left hand side by $\gamma$ ($\gamma^2$)  we see
that the left hand side is invariant under cyclic permutations
of $K,L,M$. Thus it remains to check that the factor $3$
in front of $KLM$ is correct. This however follows from
letting $K=L=M=1$.

Now fix $\beta>0$ and choose $\alpha$ and $d$ very 
small but nonzero such that
$$ 2\alpha\beta c^2+  2\beta {c}{d}=0\ \ .$$
Let
$K=1+ 2\alpha\beta {d}^2 + 2\alpha\beta d^2$,
$L=2\alpha {c} d$,
$M=2\alpha c {d}$.
Then $L^3$, $M^3$, and $3KLM$ are small of order
at least $\alpha ^2 d^2$. However,
$$K^3=1+12\alpha\beta d^2 + O(\alpha^2 d^2)\ \ \ .$$
Thus the left hand side of (\ref{counterexample})
can be made bigger than $1$.

\section{Proof of inequality (\ref{nlwc})}\label{carleson-section}

This section is very close to the known existing proofs of
Carleson's theorem in the classical linear case. We follow
closely \cite{thiele}. For example Corollary 
\ref{bessel-corollary} corresponds to a Bessel inequality
in the linear case. In the current non-linear setting it is
convenient to estimate the contribution of a single tree
pointwise outside an exceptional set (in the spirit of the
original proof by Carleson \cite{carleson}) instead of using
any $L^p$ estimate, because of the ease of pointwise
summing a geometrically decaying sequence using a quasi triangle 
inequality, see the calculation beginning with 
(\ref{geomseries}).

We shall first assume that $F$ is compactly supported.

We are interested in the dependence of constants on $d$.
It will help to introduce a constant $\Gamma$ which
(other than the constant $C$) does not change from line to line and 
has polynomial growth in $d$. The constants $C$ in this section
will be independent of $d$. 

If $G$ is a matrix in $SU(1,1)$ and $a$ is its first entry,
we shall write
$$|G|:=|a|\ \ \ .$$
Choose  $\Gamma$ so that
$$\Gamma^{-1} \log|G|\le \beta(G) \le \Gamma \log|G|\ \ .$$
By construction of $\beta$ this constant grows polynomially in $d$.

\subsection*{Orthogonality of disjoint tiles}

We define a partial ordering on tiles by $p<p'$
if $I\subset I'$ and $\omega'\subset \omega$.
Recall that all intervals $I$ and $\omega$ are
$d$- adic and assumed to be half open (containing the left but not
the right endpoint). Therefore two such intervals are
either disjoint or one is contained in the other.
Since tiles have area one we conclude that two tiles are
comparable if and only if they have non-empty 
intersection.

We observe that we have 
the following corollary of Lemma 
\ref{swapping-inequality-lemma}:

\begin{corollary}\label{bessel}

Let $\q$ be a finite set of pairwise disjoint tiles and let
$\p$ be a finite set of tiles such that for all $q\in \q$
we have 
$$q\subset \bigcup_{p\in \p}p$$
and for all $q\in \q$ and $p\in \p$ we have
$q < p$ whenever $p$ and $q$ have nonempty intersection.
Then 
\begin{equation}\label{multiswap}
\sum_{q\in \q} |I_q| \beta(G_q) \le \sum_{p\in \p} |I_p|\beta(G_{p})
\end{equation}
\end{corollary}

Proof:

If $\p$ has none or one element, then (\ref{multiswap})
is trivial because $\q$ has to be a subset of $\p$.
Fix $\p$, by induction we may assume the corollary has been 
proved for all subsets of $\p$. 
Now choose $\q$. By cancelling equal summands on both 
sides of (\ref{multiswap}) and using the result for
subsets of $\p$ we may assume that $\q$ and $\p$
are disjoint.
We may assume $\q$ is nonempty
and choose $q$ such that $l=|I_q|$ is minimal.
Since the possible values
of $l$ are discrete and bounded above by $\sum_{P\in \p}|I_p|$,
we may assume by induction that the statement
of the corollary
is true for all values of $l$
larger than a given $l_0$, and we have to prove the 
statement under the assumption  $l=l_0$. Now we
use induction on the number $n$ of tiles $q$ in 
$\q $ which satisfy $|I_q|=l$. Again
by induction we may fix an $n_0$ and assume that the
statement is true for all $n<n_0$ and we have to
prove the statement assuming $n=n_0$.

Now pick a tile $q\in \q$ such that $|I_q|=l$.
It is the vertical tile of a multitile $Q$.
We claim
\begin{enumerate}
\item $Q\subset \bigcup_{p\in \p} p$
\item Any vertical tile in $Q$ is either an element of $\q$
or it is disjoint from all tiles in $\q$.
\end{enumerate}
Assuming these two claims for now, we observe that it suffices to
prove the statement of the corollary for $\q'$ which is
the union of $\q$ and the set of vertical tiles in $Q$.
Observe that $l'=l$ and $n'<n+d-1$ where $l'$ and $n'$
are defined analogously to $l$ and $n$. By the swapping inequality
(\ref{sw}) it suffices to prove the statement for 
$\q''$ which is equal to $\q'$ with all vertical tiles of $Q$
removed and all horizontal tiles of $Q$ added in.
Observe that $l''\ge l$ and, if $l=l''$, then $n''<n$. 
Thus the statement of the corollary follows by induction.

It remains to prove the above two claims.
To see the first claim, pick $(k,x)\in Q$. There is
a $(k_0,x)\in q$. Then there is a $p\in \p$ with
$(k_0,x)\in p$. Since $p\cap q\neq \emptyset$
we have by assumptions on $\p$ that $q<p$. Thus
$I_p$ is a $d$-adic interval, strictly containing $I_q$ because of
$p\neq q$. By $d$-adicity $I_Q\subset I_p$ and hence 
$(k,x)\in p$ which had to be proved.

To see the second claim pick a vertical tile $q'$ of $Q$
and assume that $q'\cap q''\neq \emptyset$ for some $q''\in \q$. 
By minimality of the choice of $I_q$ we have
$I_{q'}\subset I_{q''}$. If this inclusion was strict,
then $q''\cap q\neq \emptyset$ which is impossible.
Hence $I_q'=I_{q''}$ and hence $q=q'$. This proves the
second claim and completes the proof of the corollary.

\endproof

\begin{corollary}\label{bessel-corollary}

Let $q$ be a set of pairwise disjoint tiles. Then
$$\sum_{q\in \q} |I_q| \beta(G_q) \le C \Gamma \int |F(x)|^2\, dx\ \ .$$
\end{corollary}

Proof: This follows from the previous corollary by a limiting
argument as in the proof of (\ref{nlwp}). 

\subsection*{Selecting trees}

We define an ordering on multitiles analogous to the ordering
on tiles. Thus $P<P'$ for two multitiles $P=I\times \omega$ and 
$P'=I'\times \omega'$
if $I\subset I'$ and
$\omega'\subset \omega$.

A set $\P$ of multitiles is called convex, if for
any three multitiles $P<P'<P''$ with $P,P''\in \P$
we can conclude $P'\in \P$.

An {\it ordered splitting} of a set $\P$ of multitiles is a 
decomposition of $\P$ into a disjoint union
$$\P=\bigcup_{n\in N}\P_n$$
where $N$ is a subset of the integers (and possibly $\infty$) 
and $P\in \P_{n}$, $P'\in \P_{n'}$ with $P<P'$ imply
$n\le n'$. Observe that if $\P$ is convex, then the
components $\P_n$ of an ordered splitting are again convex.

A {\it tree} is a set $T$ of multitiles which has a maximal element
with respect to the ordering of multitiles.
This maximal element is called the top of the tree and
denoted by $P_T$.

Each element $P$ of a tree other than the top itself has
a distinguished index $j_P\in 0,\dots, d-1$ attached to it such that
$p_{j_P}$ is the unique horizontal subtile of $P$ which 
intersects the tree top. For the top $P_T$ of a tree
we define $j_{P_T}=d-1$. Observe that we have suppressed 
the dependence of $j_P$ on the given tree in the notation.

We define the size of a collection of multitiles by
\begin{equation}\label{sizedef}
\size(\P)=\sup_{T\subset \P} |I_T|^{-1}\sum_{P\in T} \sum_{j<j_P}|I_P| \beta(G_{p_j})
\end{equation}
where the sup is taken over all trees in $\P$. For a
given tree $T$ the tiles $p_j$ occuring in the sum on the right
hand side of (\ref{sizedef}) are pairwise disjoint.

The above Corollary \ref{bessel-corollary} implies the following lemma:

\begin{lemma}\label{splitonce}

Let $\P$ be a convex set of multitiles. Then we can decompose
$\P$ into an ordered splitting $\P_1\cup\P_2$ 
such that
\begin{equation}\label{ptwo}
\size(\P_2)< 2^{-4} \size(\P)
\end{equation}
and $\P_1$ is the (not necessarily disjoint)
union of a collection $\T$ of trees such that
\begin{equation}\label{pone}
\sum_{T\in \T} |I_T| \le C \Gamma \size(\P)^{-1} \|F\|_2^2
\end{equation}
and each tree $T\in \T$ with top $P_T$ has the saturation
property that if $P<P_T$ for some $P\in \P_1$ then $P\in T$.
\end{lemma}

Proof: 
Set $\alpha=2^{-4}\size(\P)$.

We select recursively for $n=1,2,3\dots$
a tree $T_n$.
Suppose we have already chosen 
$T_m$
for all $m<n$.
If there is a tree $\tilde{T}_n$ in
$$\P^n:=\P\setminus \bigcup_{m<n} T_m$$
with size larger than $\alpha$,
then we choose one such tree with top
$P_n$ say such that the upper endpoint of $\omega_{P_n}$
is minimal. The tree $T_n$ is then the maximal
tree in $\P$ with respect to set inclusion with
top $P_n$.

We iterate this tree selection until we reach an $n=N$ such that
there is no tree in $\P^n$ with size larger than $\alpha$.
If this is the case, we stop the selection and define
$\P_1$ to be the union of trees selected. Define
$\P_2=\P\setminus \P_1$.
By maximality of each selected tree it is clear that
the splitting of $\P$ into $\P_1$ and $\P_2$ is ordered
and that each selected tree satisfies the saturation property
of the lemma.
Moreover, by the stopping condition for the tree selection 
it is clear that $\P_2$ satisfies the size estimate (\ref{ptwo}).

It remains to prove the bound (\ref{pone}). 
By Corollary \ref{bessel} it suffices to show that the set of 
tiles $p_j$ with $P\in {\tilde{T}_n}$ for some $1\le n\le N$ 
and $j<j_P$ is a set of pairwise disjoint tiles. 
Suppose to get a contradiction that $p_j<{p'}_{j'}$ for two
distinct such tiles. Then $P$ belongs to a tree $\tilde{T}_n$ and
$P'$ belongs to a tree $\tilde{T}_{n'}$. By $d$-adicity it is
easy to see that the upper endpoint of $\omega_{P_n}$ 
is greater than the upper endpoint of $\omega_{P_{n'}}$, in particular 
$n\neq n'$
and $n<n'$. But the geometry
of $p_{n'}$ qualifies it to be in the tree $T_n$, which is a 
contradiction to the maximality of $T_n$. 

This proves Lemma \ref{splitonce}.

\endproof

By iterating this lemma we obtain:

\begin{corollary}\label{splitseveral}
If $\P$ is any finite set of tiles, we can decompose it
into an ordered splitting 
$$\P=\P_{\infty}\cup \bigcup_{\k\in \Z} \P_\k$$
such that
$$\size(\P_\k)\le 2^{-4\k}$$
and $\P_\k$ is the union of a collection $\T_\k$ of trees such that
$$\sum_{T\in \T_\k}|I_T|\le C \Gamma 2^{4\k}\|F\|_2^2$$
and each tree $T\in \T_\k$ with top $P$ contains all elements
$P'\in \P_\k$ with $P'<P$.
Moreover, $\size(P_{\infty})=0$.
If the set $\P$ is convex,
then all trees in $\T_\k$ are convex.

\end{corollary}

\subsection*{A John - Nirenberg type estimate for a single tree}

Given a convex tree $T$,
we shall be concerned with the following function defined on 
$I_T$

\begin{equation}\label{maximal}
M_{T}(k)=\sup_{k\in I,I': I\subset I'\subsetneq I_T}\  
 \log \left|\prod_{P\in T:  I\subset I_P\subset I'}
\ \prod_{j<j_P} G_{p_j}(k)\right|\ \ .
\end{equation}
Here the product is to be understood in the natural
order of descending size of $I_P$ and descending $j$:
If $|I_{P}|<|I_{P'}|$, and $j=j_P>0$, $j'=j_{P'}>0$, 
then the corresponding factors appear in the order
$$\dots G_{{p'}_{j'-1}}G_{{p'}_{j'-2}}\dots G_{p'_0}\dots
G_{{p}_{j-1}}G_{{p}_{j-2}}\dots G_{p_0}\dots\ \ .$$

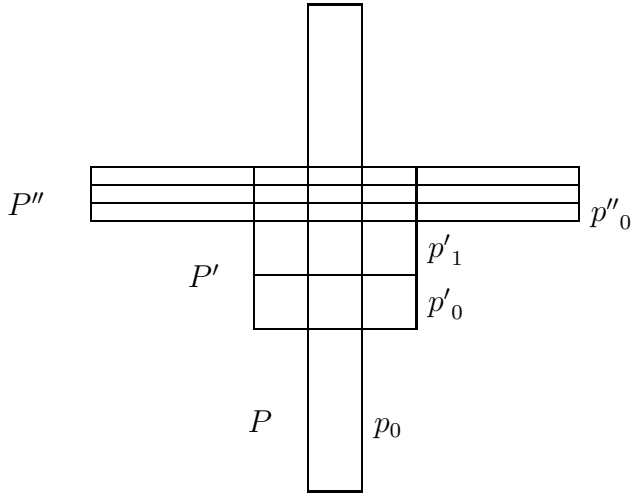
\setlength{\unitlength}{0.8mm}
\begin{picture}(160,100)

\put(70,10){\line(1,0){9}}
\put(70,37){\line(1,0){9}}
\put(70,64){\line(1,0){9}}
\put(70,91){\line(1,0){9}}

\put(70,10){\line(0,1){81}}
\put(79,10){\line(0,1){81}}

\put(60,20){$P$}
\put(81,20){$p_0$}

\put(61,37){\line(1,0){27}}
\put(61,46){\line(1,0){27}}
\put(61,55){\line(1,0){27}}
\put(61,64){\line(1,0){27}}

\put(61,37){\line(0,1){27}}
\put(88,37){\line(0,1){27}}

\put(50,44){$P'$}
\put(90,40){${p'}_0$}
\put(90,49){${p'}_1$}

\put(34,55){\line(1,0){81}}
\put(34,58){\line(1,0){81}}
\put(34,61){\line(1,0){81}}
\put(34,64){\line(1,0){81}}

\put(34,55){\line(0,1){9}}
\put(115,55){\line(0,1){9}}

\put(20,56){$P''$}
\put(117,55){${p''}_0$}

\end{picture}

We have the John-Nirenberg type lemma:

\begin{lemma}\label{jn}
Let $T$ be a convex tree. Then for every integer $\mu\ge 0$ we have
$$ |\{k\in I_T: M(k)>4 d \Gamma 2^{2\mu} \size(T)\}|\le 2^{-c \mu^2} |I_T|$$
for some small universal constant $c$.
\end{lemma}

Remark: This inequality gives less decay in $\mu$ on the right hand 
side than the usual John- Nirenberg inequality.
The loss is due to our approach to dealing with the
quasi triangle inequality in (\ref{quasi}) instead of a 
triangle inequality.

Proof:

Observe that it suffices to prove the Lemma for large $\mu$.

Let $P_T$ be the top of the tree and 
$\xi$ be the lower endpoint of the interval $\omega_{P_T}$.
Observe that all intervals $\omega_{p_j}$ for $p_j$
appearing in the product in (\ref{maximal}) lie below $\xi$. 
Therefore we do not change the value of $M_T$ if we restrict $F$ to
$[0,\xi]$. 
Therefore we shall assume for the purpose of proving this
lemma that $F$ is supported in $[0,\xi]$.

It suffices to prove a similar estimate for the simpler variant

$$\tilde{M}_T(k)=\sup_{I: k\in I} \log\left| 
\prod_{P\in T:I\subset I_P\subsetneq I_T}
\prod_{j<j_P} 
G_{p_j}(k) \right|\ \ .$$

This follows from writing
$$\prod_{P\in T:I\subset I_P\subset I'}\prod_{j<j_P} 
 G_{p_j}(k)
=
\left(\prod_{P\in T:I'\subsetneq  I_P\subsetneq I_T }\prod_{j<j_P} 
 G_{p_j}(k)\right)^{-1}
\left(\prod_{P\in T: I \subset I_P\subsetneq I_T } \prod_{j<j_P} 
 G_{p_j}(k)\right)
$$
and estimating both factors on the right hand side by $\tilde{M}_T$.
Namely, observe that for all $G,G'\in SU(1,1)$ we have
$$|G|=|G^{-1}|$$
and the quasi triangle inequality
\begin{equation}\label{quasi}
\log|GG'| \le 2\log|G| + 2\log|G'|\ \ .
\end{equation}
The latter follows from Lemma \ref{swapping}.
Thus we obtain the pointwise estimate
$${M}_T(k)\le 4\tilde{M}_T(k)\ \ ,$$ 
which reduces the matter
to estimating $\tilde{M}_T$.

We first prove for $\lambda\ge 1$ the following estimate:

\begin{equation}\label{bootstrap}
 |\{k: \tilde{M}_T(k)> \Gamma \lambda \size(T)\}|\le \lambda^{-1}  |I_T|\ \ .
\end{equation}
This follows from the estimate
\begin{equation}\label{bootstrap1}
 |\{k: \tilde{\tilde{M}}_T(k)> \lambda \size(T)\}|\le \lambda^{-1}  |I_T|\ \ .
\end{equation}
for the modified function
$$\tilde{\tilde{M}}_T(k)=\sup_{I: k\in I} \beta\left( 
\prod_{P\in T:I\subset I_P, P\neq P_T}
\prod_{j<j_P} 
G_{p_j}(k) \right)$$
because $\beta(G)$ and $\log|a|$ are comparable by a factor of $\Gamma$.

Let $P$ be a multitile of the tree $T$, let $q$ be a vertical
tile in $P$, and assume $k\in I_q$.
Then we have by support assumption
on $F$, convexity of the tree, and Lemma \ref{partition}:
$$G_{q}(k)= \prod_{P'\in T:I_P\subset I_{P'}, P'\neq P_T}\prod_{j<j_P'}
 G_{p_j}(k)\ \ . $$
Namely, it is an elementary geometric observation that the
intervals $\omega_{p_j}$ on the right form a partition
of the interval $\omega_{q}\cap (-\infty,\xi)$.

Let $\q$ be the set of maximal tiles in the set
of all tiles $q$ which are 
vertical tile of some $P\in T$
and which satisfy
$$\beta(G_{q})
\ge \lambda \size(T)\ \ .$$
Observe that the set estimated on the left hand side of
(\ref{bootstrap1}) is contained in the union of $I_q$ with $q\in \q$.
Furthermore observe that the union of all $q\in \q$ is covered by the
union of all $P\in \T$. Hence it is covered by the top multitile
$P_T$ and all horizontal tiles $p_j$ of multitiles $P\in T$ and 
$j\neq j_P$.

An application of Lemma \ref{swapping} gives
$$\sum_{q\in \q} |I_{q}|\beta(G_q)
\le \sum_{P\in \T, P\neq P_T}\ \sum_{j<j_P} |I_P|\beta(G_{p_j})\ \ .$$
Here we have used that on the right hand side we do not have to 
include the terms with $P=P_T$ or $j>j_P$ because they give
zero contribution thanks to support assumption on $F$.

This implies
$$\sum_{q\in \q} |I_q| \le \lambda^{-1} |I_T|\ \ ,$$
which proves (\ref{bootstrap1}) and therefore also 
(\ref{bootstrap}).

Now we bootstrap (\ref{bootstrap}) to the desired estimate 
for $\tilde{M}_T$. This step is analoguous to the bootstrapping
argument that can be used to prove the usual John- Nirenberg 
inequality, which is why we say Lemma \ref{jn} is of John-Nirenberg
type.
It suffices to prove:
\begin{equation}\label{step}
|\{k:\tilde M_T(k)> 2^{\mu+2} d \Gamma \size(T)\}|\le 2^{-\mu}
|\{k:\tilde M_T(k)>2^\mu d \Gamma \size(T) \}|\ \ .
\end{equation}
Namely, given this estimate, we have by iteration for every integer $\mu\ge 0$
$$|\{k:\tilde M_T(k)> 2^{2\mu}  d \Gamma \size(T)\}|\le 2^{-c \mu^2}
|\{k:\tilde M_T(k)> d \Gamma \size(T) \}|\le 2^{-c \mu^2} |I_T|\ \ .$$
This will prove the desired estimate.
We prove (\ref{step}).

Let $\P$ be the set of maximal multitiles in $P\in T$ such that
$$\max(\log|G_{q_1}|,\dots,
\log| G_{q_d}|)
\ge 2^\mu d \Gamma \size(T)$$
where $q_1,\dots,q_d$ are the vertical subtiles of $P$.
For each such $P$ let $P'\in T$ be the minimal multitile in $T$ which is larger than $P$.
(The proof trivializes if there is no such multitile because then $P=P_T$.) 
Assume that $I_P$ is equal to ${q'}_j$.

Then, by maximality,

$$\log|G_{{q'}_j}|
\le 2^{\mu} d \Gamma \size(T) $$
and moreover, for every $k\in I_P$,
$$\sup_{I: I_P\subset I } \log\left| 
\prod_{P'\in T:I\subset I_{P'}, I_{P'}\neq I_T}
\prod_{j<j_{P'}} 
G_{{p'}_j}(k) \right| \le 2^{\mu} d \Gamma \size(T)\ \ .$$

Observing that every subtree of $T$ has size at most $\size(T)$
and applying (\ref{bootstrap}) to the subtree of all $P''\in T$ with 
$I_{P''}\subset I_P$ and using the quasi triangle inequality (\ref{quasi})
we obtain: 
$$|\{k\in I_P: \tilde{M}_T(k)>2(2^{\mu}+2^{\mu}) d \Gamma \size(T)\}|\le 
2^{-\mu}d^{-1}|I_P|\ \ \ .$$
Since the intersection of $|I_P|$ and the set where 
$\tilde{M}_k>d^{\mu}\size(T)$ has measure at least $d^{-1}|I_P|$ 
(this argument is one of two places in the proof of (\ref{nlwc})
where we lose a power of $d$ in the dependence on $d$ other than 
the loss due to $\Gamma$)
we obtain:
$$|\{k\in I_T: \tilde{M}_T(k)>2^{\mu+2} d\Gamma \size(T)\}|\le 
2^{-\mu}|\{k\in I_T: \tilde{M}_k>2^{\mu} d \Gamma \size(T)\}| \ \ .
$$
This proves (\ref{step}) and completes the proof of Lemma \ref{jn}.

\subsection*{The Carleson theorem}

By restricting the set of tiles to those inside a large square 
$[0,d^K)\times [0,d^K)$ and a subsequent limiting
argument as in Section \ref{plancherel-section} we may
assume the set $\P$ of all multitiles is finite.

We are aiming to show that there exists a $C$ such that
for each $\lambda>0$ we have
$$|\{k: \sup_x \log |G(k,x)|>C d \lambda^{-1} \}|\le \Gamma^2 \|F\|_2^2
\lambda\ \ .$$
Fix $\lambda$.

Decompose $\P$ into $\P_\k$ according to Corollary \ref{splitseveral}, 
where
$$\Gamma \size(\P_\k)\le 2^{-4\k}$$
and $\P_\k$ is the union of a collection $\T_\k$ of trees
such that
$$\sum_{T\in \T_\k}|I_T|\le C \Gamma^2 2^{4\k}\|F\|_2^2$$

We define an exceptional set $E=\bigcup E_\k$.

Let $K$ be a negative integer of large modulus to be determined later.
For $\k< K$ we define
$$E_\k=\bigcup_{T\in \T_\k}I_T$$
and obtain
$$|E_\k|\le C \Gamma^2 2^{4\k}\|F\|_2^2\ \ .$$

For $\k> K$ we define
$$E_\k=\bigcup_{T\in \T_\k} \{k\in I_T:M_T(k)\ge  
4d 2^{2(\k-K)} 2^{-4\k}\}\ \ .$$
By the John-Nirenberg type Lemma  we have for $T\in \T_\k$
$$|\{k\in I_T:M_T(k)\ge 4 d 2^{2(\k-K)}  2^{-4\k}\}|
\le 2^{-c (\k-K)^2}|I_T|\ \ .
$$
Thus 
$$|E_\k|\le C \Gamma^2 2^{-c(\k-K)^2}2^{4\k}\|F\|_2^2 
= C\Gamma^2  2^{-c(\k-K)^2}2^{4(\k-K)}2^{4K}\|F\|_2^2\ \ .$$ 
If we choose $K$ maximal with 
$\lambda> C 2^{4K}$ for a certain $C$, then
$$|E|\le \sum_{\k\in \Z} |E_\k|\le \Gamma^2 \lambda\|F\|_2^2\ \ .$$

It remains to prove that for $k\notin E$ and every $x$
we have
\begin{equation}\label{outsideexc}\log |G(k,x)| \le C d \lambda^{-1}
\end{equation}
for some constant $C$.
Fix $x$.
We can write
$$G(k,x)=\prod_{P\in T} \prod_{j<j_P} G_{\omega_{p_j}}(k,x)$$
where $T$ is the convex tree of all multitiles $P$ such that
$k\in I_P$ and $x\in \omega_{P}$ and $j_P$ for $P\in T$ is the
unique index such that $x\in p_{j_P}$.

Let $T_\k$ be the intersection of $T$ with $\P_k$.
Since the sets $\P_\k$ form an ordered splitting, the sets
$T_\k$ are convex trees. 
Moreover, each $T_\k$ is contained in a tree $\tilde{T}_\k$
of $\T_\k$ by the saturation property of the trees
in $\tilde{T}_\k$.

Denote the top of $T_\k$ by $P_\k$. Observe that if $P\in T_\k$
and $P\neq P_\k$, then the number $j_P$ defined with respect to
$T$ is the same as the one defined with respect to $T_\k$.

Hence we can estimate the contribution of the tree $T_\k$
using Lemma \ref{jn} and the quasi triangle inequality
as follows 
$$
\log \left|\prod_{P\in T_\k} \prod_{j<j_P} G_{\omega_{p_j}}(k,x)\right|$$
$$\le 2 \log \left| \prod_{j<j_{P_\k}} G_{\omega_{p_j}}(k,x)\right|
+2 \log \left| \prod_{P\in T_\k, P\neq P_\k}
\prod_{j<j_{P_\k}} G_{\omega_{p_j}}(k,x)\right|$$
\begin{equation}\label{topandother}
\le 2 \log \left| \prod_{j<j_{P_\k}} G_{\omega_{p_j}}(k,x)\right|
+ 2 M_{\tilde{T}_\k}(k)\ \ .
\end{equation}
Here $p_j$ in the first summand of (\ref{topandother}) and the preceding
line
is a horizontal tile of $P_\k$.
To estimate the first term we use that $\log|G|$ is comparable
to $\beta(G)$ and the the swapping inequality to obtain
$$\log \left| \prod_{j<j_{P_\k}} G_{\omega_{p_j}}(k,x)\right|$$
$$\le \Gamma \beta(\prod_{j<j_{P_\k}} G_{\omega_{p_j}}(k,x))$$
$$\le d \Gamma \sum_{j<j_{P_\k}} \beta(G_{\omega_{p_j}})\le d 2^{-4\k}$$
The last inequality follows by observing that $\{P_\k\}$ constitutes a tree by itself
which is controlled in size because $P\in \P_\k$.
Observe that in this argument we lose a factor $d$.

Thus, by choice of $x$.
$$\log \left|\prod_{P\in T_\k} \prod_{j<j_P} G_{\omega_{p_j}}(k,x)
\right|
\le C d\Gamma 2^{2(\k-K)} 2^{-4\k}\\ \ .$$

The trees $T_\k$ with $\k<K$ are empty, because $x$ is not in the exceptional
set. 
Moreover, by finiteness assumption on the set $\P$
there is a $K'$ so that $T_\k$ is empty for $\k>K'$ and $\k\neq \infty$. 
(The constant $C$ is not allowed to depend on $K'$)

Then we have
$$G(k,x)=\left(
\prod_{P\in T_\infty}
\prod_{j<j_P} G_{\omega_{p_j}}(k,x)\right)
\left(\prod_{K\le \k\le K'} \prod_{P\in T_\k}
\prod_{j<j_P} G_{\omega_{p_j}}(k,x)\right)
$$
where as usual the product has to be read in the correct order.
The factor coming from $T_\infty$ can be discarded since it gives a 
unitary matrix.

By the quasi triangle inequality we have
$$\log|G(k,x)|\le 4 C d 2^{-4K}\ \ .$$
Namely, we can prove inductively
\begin{equation}\label{geomseries}
\log|\prod_{K''\le \k\le K'} \prod_{P\in T_\k}
\prod_{j<j_P} G_{\omega_{p_j}}(k,x)
|
\end{equation}
$$\le 2\log|\prod_{P\in T_{K''}}
\prod_{j<j_P} G_{\omega_{p_j}}(k,x)
|
+2\log |\prod_{K''+1\le \k\le K'} \prod_{P\in T_\k}
\prod_{j<j_P} G_{\omega_{p_j}}(k,x)
|$$
$$\le 2 C d 2^{2(K''-K)} 2^{-4K''}
+  2 C d 2^{2(K''-K)} 2^{-4K''}
$$
$$\le C d 2^{2((K''-1)-K))} 2^{-4(K''-1)}\ \ ,$$
where contrary to our standing convention the constant
$C$ for induction purpose is the same in all 
appearances in this calculation.

This proves inequality (\ref{outsideexc}) and therefore
completes the proof of inequality (\ref{nlwc})
in the case of compactly supported $F$.

It remains to discuss the case of not necessarily compactly
supported  $F\in L^2(\R)$. 
We first prove that the limit (\ref{limit}) exists
almost everywhere. 
It suffices to fix small $\epsilon$
and prove that the limit exists for all $k$ outside a set
of measure $\epsilon$. 
This can be done by decomposing the positive real axis into
intervals $[0,x_1)$, $[x_1,x_2)$, etc, such that
the $L^2$ norms of the restrictions of $F$ to the intervals
$\omega_j=[x_j,x_{j+1})$
decay very rapidly in $j$. This implies by the Plancherel
inequality (\ref{nlwp}) that for 
$k$ outside a small set (of size $\epsilon/2$) the values
$\log|G_{\omega_j}|$ are still very rapidly decaying, so that
one can use the triangle inequality 
\begin{equation}\label{idealtriangle}
\log|G_1G_2|\le \log|G_1| +  \log \|G_2\|_{op}
\end{equation}
to show that 
the sequence $\log|G_{[0,x_j)}|$ is a Cauchy sequence.
Using (\ref{nlwc}) for for the restriction of $F$ to
each interval $[x_j,x_{j+1})$ one can observe that the
operator norms of the matrices $G_{[x_j,x)}$ with $x_j<x<x_{j+1}$
are small for large $j$ and all $k$ outside a set of measure
$\epsilon/2$. Using the triangle 
inequality (\ref{idealtriangle}) one can show that for $k$ outside
a set of measure $\epsilon$ the limit of $\log|G_{[0,x)}|$ exists. 
This proves existence of the limits in (\ref{limit}). Similar arguments 
as these make it straight forward to prove (\ref{nlwc}) for arbitrary 
potentials $F\in L^2([0,\infty))$.

\section{Multilinear expansions}\label{multilinear-section}

Writing the differential equation
(\ref{wode}) as an integral equation
$$G(k,x)= 
\left(\begin{array}{cc} 1  & 0 
\\ 0 & 1 \end{array}\right)
+
\int_{-\infty}^x
\left(\begin{array}{cc} 0  & F(x) w(k,t)
\\ \overline{F(x)} \overline{w(k,t)} & 0 \end{array}\right)
G(k,t)\, dt
$$
we can use Picard iteration to obtain the formal
solution
\begin{equation}\label{walshpicard}
G(k,x)=
\left(\begin{array}{cc} 1  & 0 
\\ 0 & 1 \end{array}\right)
+
\sum_{n=1}^\infty \int_{t_1<\dots<t_n<x}\prod_{j=n}^1 
\left(\begin{array}{cc} 0  & F(t_j) w(k,t) 
\\ \overline{F(t_j)} \overline{w(k,t)} 
& 0 \end{array}\right)\, dt_j\ \ .
\end{equation}
Christ and Kiselev \cite{ck0},\cite{ck1}
prove convergence for almost every $k$ of this formal expansion
if $F\in L^p(\R)$ with $p<2$, and they use the expansion to show 
the maximal Hausdorff Young inequality (they work on a different
model of the nonlinear Fourier transform, but their arguments apply
to this case too). In \cite{coex} it has been shown that the higher order terms of
the Fourier analogue of this expansion are unbounded for $F\in L^2$ and
therefore not very well suited to be used to prove a nonlinear 
Carleson theorem.

In this section we show that a similar discussion as in \cite{coex}
applies in the $d$-adic setting provided $d\ge 3$.
More precisely we will focus on the quadratic term in the above expansion
and on the case $d=3$ and prove Proposition \ref{unbounded-3} below.
The arguments generalize to $d>3$ and to the higher-linear terms,
but we shall not elaborate on this because our main point can be made
clear for this special case.

The quadratic term in the above expansion
is a diagonal matrix with entries
\begin{equation}\label{quadrterm}
Q(F)(k,x)=\int_{t_1<t_2<x} 
{F(t_1)} {w(k,t_1)}
\overline{F(t_2)  w(k,t_2)}
\, dt_1dt_2
\end{equation}
and the complex conjugate of (\ref{quadrterm}).

We consider 
\begin{equation}\label{maxquadr}
M(F)(k)=\sup_x \left|
\int_{t_1<t_2<x} F(t_1) 
w(k,t_1) 
\overline{F(t_2) w(k,t_2)}\, dt_1dt_2\right|\ \ .
\end{equation}

The following proposition implies that there is no reasonable
a priori bound for the size of this function
in terms of the $L^2$ norm
of  $F$.

\begin{proposition}\label{unbounded-3}
Let $d=3$.
There is an $\epsilon>0$ such that for each 
$N>0$ there is a finite $d$-adic step function $F$ with $L^2$
norm $1$ such that
$$\|\{k: M(F)(k)>N\}\|>\epsilon\ \ \ .$$
\end{proposition}

Proof:

Pick a large integer $N$. In the proof of Lemma
\ref{constance} we have seen that for a given tile
$p=I\times \omega$ we can write
$$w(k,x)1_I(k) 1_\omega(x)= w_p(k)\overline{\widehat{w}_p(x)}$$
for some functions $w_p$ and $\widehat{w}_p$ which have constant
modulus on $I$ and $\omega$ respectively and which we may 
assume to have $L^2$ norm $1$. We observe for $x\in \omega$:
\begin{equation}\label{packettransform}\int w_p(k) \overline{w(k,x)} \, dk= 
\int w_p(k) \overline{w_p(k)} \widehat{w}_p(x) \, dk= \widehat{w}_p(x)\ \ .
\end{equation}
Thus, relying on the well known fact that the Cantor group 
Fourier transform
is an isometry in $L^2$ (this can be shown by a linearized version of the
arguments in Section \ref{plancherel-section}), we see that the integral
on the left hand side of (\ref{packettransform}),
which is the Walsh-Fourier transform of $w_p$, has to vanish
outside $\omega_P$ and thus $\widehat{w_p}$ 
is indeed the Cantor group Fourier transform of $w_p$.

Consider the tiles
$$p_j= I_j\times \omega_j:=[3^{-N}j,3^{-N} (j+1))\times [3^{N}j,3^{N} (j+1))$$
for $j=0,\dots, 3^N-1$ and set
$$F(x):=\sum_{j=0}^{3^N-1}F_j(x):=\sum_{j=0}^{3^N-1} 3^{-N/2} \widehat{w}_{p_j}(x)\ \ .$$
Thus $\|F\|_2=1$.

Observe that the intervals $\omega_j$ above form a partition
of the interval $[0,3^{2N})$. For $k\in [0,1]$ let $x(k)$ be the left 
endpoint of the
unique interval $\omega_j$ which contains $3^{2N}k$.
Consider for $k\in [0,1)$:
$$\tilde{Q}(F)(k)=\int_{t_1<t_2<x(k)} 
F(t_1)  w(k,t_1)\overline{F(t_2)w(k,t_2)}\, dt_1dt_2$$
We shall show that the imaginary part of $\tilde{Q}$ is large
on a big set:
$$\{k\in [0,1): |\Im(\tilde{Q}(F)(k))|> N/4\}>  1/3\ \ ,
$$
which will prove the proposition. Our argument works only for the 
imaginary part. Indeed, the real part of $\tilde{Q}$ can be seen
to satisfy good estimates. This is the reason why our argument works
only for $d\ge 3$. If $d=2$, then the characters $w(k,x)$ are 
purely real and our argument does not work. We do not know
whether for $d=2$ the series (\ref{walshpicard}) converges 
for genuinly complex $F\in L^2$ in a reasonable sense.

We may split
$$\tilde{Q}(F)(k)=\sum_{0\le j,j'\le 3^{N}-1}\int_{t_1<t_2<x(k)} 
F_j(t_1)  w(k,t_1)\overline{F_{j'}(t_2)w(k,t_2)}\, dt_1dt_2$$
If $j>j'$ then the integrand is zero on the domain
of integration, thus we may disregard these terms.
Likewise, if $3^Nj'\ge x(k)$, the integrand is zero.
If $j<j'< 3^{-N} x(k)$, then the constraint $t_1<t_2<x(k)$ in the domain
of integration is superfluous and we can write
$$\int_{t_1<t_2<x(k)} 
F_j(t_1) w(k,t_1)\overline{F_{j'}(t_2)w(k,t_2)}\, dt_1dt_2$$
$$=\int  F_j(t_1) w(k,t_1) \, dt_1 \ 
\int \overline{F_{j'}(t_2) w(k,t_2)}\, dt_2$$
$$=w_{p_j}(k)\overline{w_{p_{j'}}(k)}=0\ \ .$$
Thus it only remains to consider the terms with
$j=j'<3^{-N}x(k)$:
$$\tilde{Q}(F)(k)=\sum_{0\le j< 3^{-N}x(k)}\int_{t_1<t_2} 
F_j(t_1) w(k,t_1)\overline{F_j(t_2)w(k,t_2)}\, dt_1dt_2\ \ .$$

For each $ t_1<t_2 $, $t_1,t_2\in \omega_j$
there is a minimal $3$-adic interval $\omega\subset \omega_j$ such that
$t_1$ and $t_2$ are both contained in $\omega$.
Then $t_1$ and $t_2$ are in different $3$-adic subintervals $\omega_{(m)}$
and $\omega_{(m')}$ of the next smaller generation of $\omega$ . 
Indeed, $m<m'$.
We split $\tilde{Q}(F)$ according to the size of $\omega$ as follows:
$$\tilde{Q}(F)(k)=\sum_{\kappa\le N}\ 
\sum_{0\le j< 3^{-N}x(k)}\ \sum_{|\omega|=3^\kappa}\ \sum_{0\le m<m'\le 2}
\int_{\omega_{(m)}}  F_j(t_2) w(k,t_2) \, dt_2 \ 
\int_{\omega_{(m')}} \overline{F_j}(t_1) \overline{w(k,t_1)}\, dt_1$$

Now fix $\kappa$ and $j$.
Let the $-\kappa$-th coefficient in the ternary
expansion of $k$ be $k_{-\kappa}$ and the $-\kappa$-th coefficient in
the ternary expansion of $3^{-N} j$ be $j_{N-\kappa}$. 
Then we observe that for $m=0,1,2$
$$\int_{\omega_{(m)}}  F_j(t_2) w(k,t_2) \, dt_2 =
\gamma^{m k_\kappa-mj_{N-\kappa}}\int_{\omega_{(0)}}  
F_j(t_2) w(k,t_2) \, dt_2$$
since there is a bijection of $\omega_{(0)}$ to $\omega_{(m)}$
given by switching the $\kappa-1$-st coefficient of each element
in $\omega_{(0)}$ from $0$ to $m$. Moreover
$$ |\int_{\omega_{(m)}}  F_j(t_2) w(k,t_2) \, dt_2| = 3^{-N+\kappa-1}$$
provided $\omega \subset [3^Nj,3^N(j+1))$ and $3^{-N}j$ and $k$
are in the same $3$-adic interval $I_{(j)}$ of length $|\omega|^{-1}$. 
Otherwise the integral on the left hand side is zero. Thus
\begin{equation}\label{single-scale}
\Im \sum_{|\omega|=3^\kappa}\ \sum_{1\le m<m'\le 3}
\int_{\omega_{(m)}}  F_j(t_2) w(k,t_2) \, dt_2 \ 
\int_{\omega_{(m')}} \overline{F_j}(t_1) \overline{w(k,t_1)}\, dt_1
\end{equation}
$$=1_{I_{(j)}}(k)3^{-N+\kappa-2}\Im \sum_{1\le m<m'\le 3}\gamma^{(m-m')
(k_\kappa-j_{N-\kappa})}
$$
The imaginary part of the sum 
\begin{equation}\label{m-sum}\sum_{1\le m<m'\le 3}\gamma^{(m-m')
(k_\kappa-j_{N-\kappa})}
\end{equation}
is equal to $0$ if $k_\kappa-j_{N-\kappa}=0$, it is
equal to $-\Im(\gamma)$ if $k_\kappa-j_{N-\kappa}=1$ or $-2$, 
and it is equal to $\Im(\gamma)$ if $k_\kappa-j_{N-\kappa}=-1$ or $2$.

Now we add the terms (\ref{single-scale}) over all $0<j<3^{-N}k$,
i.e. we consider
\begin{equation}\label{all-j}
\Im \sum_{0\le j< 3^{-N}x(k)}\ \sum_{|\omega|=3^\kappa}\ \sum_{0\le i<i'\le 2}
\int_{\omega_{i}}  F_j(t_2) w(k,t_2) \, dt_2 \ 
\int_{\omega_{i'}} \overline{F_j}(t_1) \overline{w(k,t_1)}\, dt_1\ \ .
\end{equation}

By the previous remarks it suffices to count the terms for which the 
contribution (\ref{single-scale}) is equal to $-\Im(\gamma)$
and the number of terms for which it is equal to $\Im(\gamma)$.
If $k$ is in the left third of a $3$-adic interval
of length $3^{-\kappa}$ then we do not have a nonzero
contribution from any $j$ because the constraint $3^{-N}j< k$ and the
constraint that $3^{-N}j$ and $k$ are in the same $3$-adic interval
of length $3^{-\kappa}$ implies that $3^{-N}j$ is in the left
third of the same $3$-adic interval of length $3^{-\kappa}$ as $k$
and thus $j_{N-\kappa}=0$. By the discussion of (\ref{m-sum}) 
we therefore see that (\ref{all-j}) is equal to $0$.

If $k$ is in the middle third of a $3$-adic interval of length $3^{-\kappa}$,
then we get a contribution if $3^{-N}j$ is in the left third of that interval.
There are $3^{N-\kappa-1}$ such values of $j$, thus (\ref{all-j}) is equal
to $-3^{-3}\Im(\gamma)$.

If $k$ is in the right third of a $3$-adic interval of length $3^{-\kappa}$, 
then we get as many
$j$ in the left third as in the middle third, and their
contributions cancel each other, thus (\ref{all-j}) is
equal to $0$.

Now summing (\ref{all-j}) over all $\kappa\le N$ reduces to
counting the number of scales $\kappa$ for which $k$ is in the
middle third of a $3$-adic interval of length $3^{-\kappa}$.
We may restrict attention to $\kappa\ge 0$ because we assume
$k\in [0,1)$ and hence $k$ is always in the left third
of any $3$-adic interval of length $3^{-\kappa}$ if $\kappa<0$.

Thus
$$\Im(\tilde{Q}(F))(k)=-3^{-3}\Im(\gamma) \# \{0\le \kappa \le N: 
k_{-\kappa-1}=1\}\ \ .$$
Since for each scale exactly one third (in measure) of all numbers in $[0,1)$
are in the middle interval of a $3$-adic interval of that scale, we obtain
$$\int_0^1 \Im(\tilde{Q}(F))(k)=-(N+1) 3^{-4}\Im(\gamma)\ \ .$$
Moreover, clearly
$$\sup_{k\in [0,1)}|\Im(\tilde{Q}(F))(k)|=(N+1)3^{-3}\Im(\gamma)\ \ .$$
Thus
$$\|\{k: |\Im(\tilde{Q}(F))(k)|\ge (N+1)3^{-5}\Im(\gamma)\}\|\ge 3^{-2}\ \ .$$
This proves \ref{unbounded-3} since the choice of $N$ was arbitrary.
\endproof

\section{Appendix}\label{appendix-section}

In this section we prove the formula (\ref{fourier-plancherel}).
The proof uses complex contour integration, which is essentially
the only method we know to prove the inequality. It is a global 
argument, which should be contrasted to the argument 
in Section \ref{plancherel-section} for the Cantor group case,
which uses local methods. The local methods are useful in
proving Carleson's theorem, while global methods are hard to adapt.

The argument in this section is well known, variants of it go
back at least as far as the article by Buslaev and Faddeev 
\cite{buslaevfaddeev} or the work by Verblunsky 
\cite{verblunsky1},\cite{verblunsky2} in the discrete case.

Let $F$ be a compactly supported, complex valued,
smooth function on $\R$.
Consider the solution to (\ref{equivalent-ode})
with initial condition $G(-\infty)=\id$.
Writing this initial value problem as an integral equation and using Picard iteration
as in Section \ref{multilinear-section} gives the solution as
a formal expansion
$$G(k,x)=
\left(\begin{array}{cc} 1  & 0 
\\ 0 & 1 \end{array}\right)
+
\sum_{n=1}^\infty \int_{t_1<\dots<t_n<x}\prod_{j=n}^1 
\left(\begin{array}{cc} 0  & F(t_j) e^{2 ikt_j} 
\\ \overline{F(t_j)} e^{-2 i kt_j} & 0 \end{array}\right)\, dt_j\ \ .
$$
Indeed, this expansion is easily seen to converge
using that the $L^1$ norm of $F$ is finite and a symmetry
argument for the integration domain to obtain a factor $1/n!$
for the $n$-th multilinear term.
At $x=\infty$ we obtain for the first entry $a(k)$
of $G(k,\infty)$
$$a(k)=1+\sum_{n=1}^\infty \int_{t_1<t_2<\dots<t_{2n}} 
\prod_{j=1}^n F(t_{2j-1})\overline{F(t_{2j})}e^{2 i k (t_{2j}-t_{2j-1})}\, dt_{2j-1}
dt_{2j}\ \ .$$
This function $a(k)$ extends holomorphically to $k$ in the half plane $\Im(k)\ge 0$ because it is a summable superposition of functions of the form $e^{ikt}$
with $t>0$.
Assume for now that $a(k)$ does not have any zeros in the closed upper half
plane, we will prove this at the end of this section.

Then we can define a function $\log(a(k))$ in the upper half plane.
We choose the branch of the logarithm so that for $|k|\to \infty$
we have $\log(a(k))\to 0$. It will become clear momentarily
that this is well defined.

We consider the counter clockwise contour integral over a large semicircle
$C=C_1+C_2$ where $C_1=[-r,r]$ and $C_2=\{k: |k|=r, Im(k)\ge 0\}$.
We show that on $C_2$ only the
first nontrivial term in the expansion of $a(k)$ gives a contribution
to the integral.
We do a partial integration for this term
$$
\int_{t_1<t_2} 
F(t_{2})\overline{F(t_{1})}e^{2 i k(t_{2}-t_{1})}\, dt_{2}dt_1$$
$$
=\int_{s>0} \int_t
F(t+s)\overline{F(t)}e^{2 i ks}\, dtds$$
$$
=-\frac 1{2 i k} \int_t
F(t)\overline{F(t)}\, dt
-\frac 1{2 i k} \int_{s>0}\int_{t}
F'(t+s)\overline{F(t)}e^{2 i ks}\, dtds$$
$$=-\frac 1{2 i k} \|F\|_2^2 +O(k^{-2})\ \ .$$
Here the estimate on the remainder term can be seen by one further
partial integration.

For all the other terms in the expansion 
we can do partial integration in all variables $s_j=t_{2j-1}-t_{2j}$
so as to get the estimate $O(k^{-2})$ or even better for all these terms.
Thus we have for large $|k|$
$$\log(a(k))= - \frac 1{2 ik}\|F\|_2^2 +O(k^{-2})\ \ .$$
Doing the integration on $C_2$ we obtain
$$\int_{C_2} \log(a(k))\, dk= -\frac \pi 2 \|F\|_2^2+O(r^{-1})\ \ .$$
Since the contour integral over $C$ vanishes, we obtain in the limit $r\to \infty$
$$\int_{\R} \log(a(k))\, dk= \frac \pi 2 \|F\|_2^2\ \ ,$$
which implies (\ref{fourier-plancherel}) because the right hand side is real.

It remains to prove that $a$ does not have any
zeros in the upper half plane. To this end write 
$a(k,x)$ and $b(k,x)$ for the entries in the first
column of $G(k,x)$ and consider
the quantity
\begin{equation}\label{monotone}
|a(k,x)|^2|e^{-2ikx}|^2 - |b(k,x)|^2\ \ .
\end{equation}
Writing $\Re(k)=\sigma$ and $\Im(\kappa)=\tau$,
the partial derivative in the $x$ variable 
of this expression is
$$2\Re\left[F(x) e^{2i\sigma x} e^{2 \tau x} b(k, x) \overline{a}(k,x)
- \overline{F}(x) e^{-2i\sigma x} e^{2 \tau x} a(k, x) \overline{b}(k,x)\right]
+ 2 \tau |a(k,x)|^2 e^{2 \tau x}$$
$$= 4 \tau |a(k,x)|^2 e^{4\tau x}\ \ .$$
The latter is always positive for $\tau>0$. Since (\ref{monotone})
is equal to 
$$|e^{-2ikx}|^2$$
for $x$ near $-\infty$, we conclude that (\ref{monotone}) is
positive for all $x$ and all $\tau>0$.
This proves that $a(x,k)$ is nonzero for such $x$ and $\tau$.
For $\tau=0$ we observe by a similar argument
that $|a(k,x)|^2-b(k,x)|^2$ is constant equal to $1$ and thus $a$
has no zeros on the real axis neither.

\end{document}